\documentclass[a4paper,12pt]{article}
\usepackage{amsmath, amsfonts, amssymb}
\usepackage{xcolor}
\usepackage[all]{xy}

\makeatletter
\newbox\bk@bxb
\newbox\bk@bxa
\newif\if@bkcont
\newcount\bk@lcnt

\def\breakboxskip{2pt}
\def\breakboxparindent{1.8em}

\def\breakbox{\vskip\breakboxskip\relax
\setbox\bk@bxb\vbox\bgroup
\advance\linewidth -2\fboxrule
\hsize\linewidth\@parboxrestore
\parindent\breakboxparindent\relax}

\def\bk@split{%
\@tempdimb\ht\bk@bxb 
\advance\@tempdimb\dp\bk@bxb
\setbox\bk@bxa\vsplit\bk@bxb to\z@ 
\setbox\bk@bxa\vbox{\unvbox\bk@bxa}
\setbox\@tempboxa\vbox{\copy\bk@bxa\copy\bk@bxb}
\advance\@tempdimb-\ht\@tempboxa
\advance\@tempdimb-\dp\@tempboxa}

\def\bk@addfsepht{%
\setbox\bk@bxa\vbox{\vskip\fboxsep\box\bk@bxa}}

\def\bk@addskipht{%
\setbox\bk@bxa\vbox{\vskip\@tempdimb\box\bk@bxa}}

\def\bk@addfsepdp{%
\@tempdima\dp\bk@bxa
\advance\@tempdima\fboxsep
\dp\bk@bxa\@tempdima}

\def\bk@addskipdp{%
\@tempdima\dp\bk@bxa
\advance\@tempdima\@tempdimb
\dp\bk@bxa\@tempdima}

\def\bk@line{%
\hbox to \linewidth{%
\hskip-2\fboxsep\vrule \@width\fboxrule\hskip.5\fboxsep\vrule \@width\fboxrule\hskip1.5\fboxsep
\box\bk@bxa\hfil
}}%

\def\endbreakbox{\egroup
\ifhmode\par\fi{\noindent\bk@lcnt\@ne
\@bkconttrue\baselineskip\z@\lineskiplimit\z@
\lineskip\z@\vfuzz\maxdimen
\bk@split\bk@addfsepht\bk@addskipdp
\ifvoid\bk@bxb 
\def\bk@fstln{\bk@addfsepdp
\hskip-\parindent\vbox{\llap{\raisebox{-2ex}{\rule{1.5\fboxsep}{\fboxrule}\hskip.5\fboxsep}}\bk@line\llap{\rule{1.5\fboxsep}{\fboxrule}\hskip.5\fboxsep}}}

\else 
\def\bk@fstln{\vbox{\llap{\raisebox{-2ex}{\rule{1.5\fboxsep}{\fboxrule}\hskip.5\fboxsep}}\bk@line}\hfil%
\advance\bk@lcnt\@ne
\loop
\bk@split\bk@addskipdp\leavevmode
\ifvoid\bk@bxb 
\@bkcontfalse\bk@addfsepdp
\vtop{\bk@line\llap{\rule{2\fboxsep}{\fboxrule}}}%

\else 
\bk@line
\fi
\hfil\advance\bk@lcnt\@ne
\if@bkcont\repeat}%
\fi
\leavevmode\bk@fstln\par}\vskip\breakboxskip\relax}

\def\smp{\smallskip\par}
\def\un{{\bf 1}}
\def\zero{\{0\}}
\def\pf{\noindent{\bf Proof:}\ }

\def\findemo{~\leaders\hbox to 1em{\hss\  \hss}\hfill~\raisebox{.5ex}{\framebox[1ex]{}}\smp}

\def\mpn{\medskip\par\noindent}
\def\smpn{\smallskip\par\noindent}
\def\normal{\mathop{\trianglelefteq}}

\def\smp{\smallskip\par}
\def\smpn{\smallskip\par\noindent}

\def\mpoint{\;\;.}
\def\mvirg{\;\;,}

\def\Res{{\rm Res}}
\def\Ind{{\rm Ind}}
\def\Inf{{\rm Inf}}
\def\Def{{\rm Def}}
\def\Iso{{\rm Iso}}

\def\Hom{{\rm Hom}}

\def\Inf{{\rm Inf}}

\def\Ker{{\rm Ker}}

\def\Id{{\rm Id}}

\def\op{^{op}}

\def\Z{\mathbb{Z}}

\def\F{\mathbb{F}}

\newcommand{\romain}[1]{\uppercase\expandafter{\romannumeral #1}}

\newcommand{\flh}[2]{\mathop{\hbox to 12mm{\rightarrowfill}}_{\displaystyle #2}^{\displaystyle #1}\limits}
\newcommand{\sflh}[2]{\mathop{\hbox to 12mm{\rightarrowfill}}_{\scriptstyle #2}^{\scriptstyle #1}\limits}

\newcommand{\sou}[1]{\underline{#1}}
\newcommand{\sur}[1]{\,\overline{\! #1}}

\def\op{^{op}}

\newcommand{\carre}[8]{\begin{array}{ccc}
#1&\mathop{\hbox to 12mm{\rightarrowfill}}^{\displaystyle{#2}}\limits&#3\\
\llap{$\displaystyle{#4}$}\left\downarrow\vbox to 6mm{}\right. & & \left\downarrow\vbox to 6mm{}\right.\rlap{$\displaystyle{#5}$}\\
#6&\mathop{\hbox to 12mm{\rightarrowfill}}_{\displaystyle #7}\limits&#8\\
\end{array}}

\newcommand{\carrem}[8]{\begin{array}{ccc}
#1&\mathop{\hbox to 12mm{\rightarrowfill}}^{\displaystyle #2}\limits&#3\\
\llap{$\displaystyle #4$}\left\uparrow\vbox to 6mm{}\right. & & \left\uparrow\vbox to 6mm{}\right.\rlap{$\displaystyle #5$}\\
#6&\mathop{\hbox to 12mm{\rightarrowfill}}_{\displaystyle #7}\limits&#8\\
\end{array}}

\newenvironment{enonce}[1]{\pagebreak[2]\refstepcounter{subsection}\refstepcounter{prop}\smpn{{\bf \thesection.\arabic{prop}.\ \ #1:}}\begin{it} }{\end{it}\smp}
\newenvironment{enonce*}[1]{\pagebreak[2]\smpn{#1:}\begin{it} }{\end{it}\smp}
\newcommand{\result}[1]{\begin{enonce}{#1}}
\def\fresult{\end{enonce}}

\newenvironment{mth}[1]{\begin{breakbox}\begin{enonce}{#1}}{\end{enonce}\end{breakbox}}
\newenvironment{mth*}[1]{\begin{breakbox}\begin{enonce*}{#1}}{\end{enonce*}\end{breakbox}}
\newenvironment{rem}[1]{\refstepcounter{subsection}\refstepcounter{prop} \mpn{{\bf \thesection.\arabic{prop}.}\ \ \bf#1:}}{\smp}

\makeatletter
\renewenvironment{enumerate}{\ifnum \@enumdepth >3 \@toodeep\else
      \advance\@enumdepth \@ne
      \edef\@enumctr{enum\romannumeral\the\@enumdepth}\list
      {\csname label\@enumctr\endcsname}{\setlength{\topsep}{1ex}\setlength{\itemsep}{0pt}\usecounter
        {\@enumctr}\def\makelabel##1{\hss\llap{##1}}}\fi}{\endlist}
\renewenvironment{itemize}{\ifnum \@itemdepth >3 \@toodeep\else \advance\@itemdepth \@ne
\edef\@itemitem{labelitem\romannumeral\the\@itemdepth}%
\list{\csname\@itemitem\endcsname}{\setlength{\topsep}{1ex}\setlength{\itemsep}{0pt}\def\makelabel##1{\hss\llap{##1}}}\fi}
{\endlist}
\makeatother
\makeatletter
\def\@sect#1#2#3#4#5#6[#7]#8{\ifnum #2>\c@secnumdepth
    \let\@svsec\@empty\else
    \refstepcounter{#1}\edef\@svsec{\csname the#1\endcsname .\hskip .5em}\fi
    \@tempskipa #5\relax
     \ifdim \@tempskipa>\z@
       \begingroup #6\relax
         \@hangfrom{\hskip #3\relax\@svsec}{\interlinepenalty \@M #8\par}%
       \endgroup
      \csname #1mark\endcsname{#7}\addcontentsline
        {toc}{#1}{\ifnum #2>\c@secnumdepth \else
                     \protect\numberline{\csname the#1\endcsname}\fi
                   #7}\else
       \def\@svsechd{#6\hskip #3\relax  
                  \@svsec #8\csname #1mark\endcsname
                     {#7}\addcontentsline
                          {toc}{#1}{\ifnum #2>\c@secnumdepth \else
                            \protect\numberline{\csname the#1\endcsname}\fi
                      #7}}\fi
    \@xsect{#5}}
\def\section{\@startsection {section}{1}{\z@}{-3.5ex plus-1ex minus
    -.2ex}{2.3ex plus.2ex}{\reset@font\Large\bf}}  

\makeatother
\renewenvironment{equation}{\refstepcounter{subsection}\refstepcounter{prop}$$}{\leqno{\bf (\theprop)}$$}

\def\mar[#1]{\ar@{-}[#1]|-{\object@{<}}}
\def\marb[#1]{\ar@{-}[#1]|{\object+{  }}}

\def\F{\mathbb{F}}
\def\FB{\mathbb{F}B}

\def\CP{\mathcal{P}}
\def\endpf{\findemo}
\newcommand{\BgrK}[1]{\mathcal{B}_{#1}{\rm\hbox{-}gr}}
\newcommand{\BgrKp}[2]{\mathcal{B}_{#1}^{(#2)}{\rm\hbox{-}gr}}
\newcommand{\grpover}[1]{\mathsf{grp}_{\Downarrow#1}}
\newcommand{\AMod}[1]{#1\hbox{-}\mathrm{Mod}}

\begin{document}
\centerline{\Large\bf Relative $B$-groups}\vspace{.5cm}\par
\centerline{\bf Serge Bouc}\vspace{1cm}\par
{\footnotesize {\bf Abstract :}  This paper extends the notion of $B$-group to a relative context. For a finite group $K$ and a field $\F$ of characteristic 0, the lattice of ideals of the Green biset functor $\FB_K$ obtained by shifting the Burnside functor $\FB$ by $K$ is described in terms of {\em $B_K$-groups}. It is shown that any finite group $(L,\varphi)$ over $K$ admits a {\em largest quotient $B_K$-group} $\beta_K(L,\varphi)$. The simple subquotients of $\FB_K$ are parametrized by $B_K$-groups, and their evaluations can be precisely determined. Finally, when $p$ is a prime, the restriction $\FB_K^{(p)}$ of $\FB_K$ to finite $p$-groups is considered, and the structure of the lattice of ideals of the Green functor $\FB_K^{(p)}$ is described in full detail. In particular, it is shown that this lattice is always finite.
\vspace{2ex}\par}
{\footnotesize {\bf AMS Subject classification :} 18B99, 19A22, 20J15.\vspace{2ex}}\par

\section{Introduction}
In the study of the lattice of biset-subfunctors of the Burnside functor $\FB$ over a field $\mathbb{F}$ of characteristic 0 (cf. Section~7.2 of \cite{doublact}, or Chapter~5 of~\cite{bisetfunctorsMSC}), a special class of finite groups, called {\em $B$-groups}, plays an important role: indeed, the simple subquotients of the biset functor $\mathbb{F}B$ are exactly the functors $S_{H,\mathbb{F}}$, where $H$ is such a $B$-group. It was shown moreover that each finite group $G$ has a largest quotient $B$-group $\beta(G)$. \par
Let $K$ be a fixed finite group. This paper proposes a generalization of the above methods and notions, in order to study the lattice of ideals of the {\em shifted Burnside functor} $\FB_K$. We start by introducing a category $\grpover{K}$ of groups over $K$, similar to the comma category of finite groups over $K$, in which morphisms are obtained by allowing diagrams to commute up to inner automorphisms of $K$. \par
To each such group $(L,\varphi)$, where $\varphi:L\to K$, is attached a specific ideal $\mathsf{e}_{L,\varphi}$ of $\FB_K$, and it is shown that every ideal of $\FB_K$ is equal to the sum of the ideals $\mathsf{e}_{L,\varphi}$ it contains. A special class of groups over $K$ is introduced, called $B_K$-groups, and it is shown that for each group $(L,\varphi)$ over $K$, there exists a largest $B_K$-group $\beta_K(L,\varphi)$ quotient of $(L,\varphi)$. Moreover $\mathsf{e}_{L,\varphi}=\mathsf{e}_{\beta_K(L,\varphi)}$. It follows that the lattice of ideals of $\FB_K$ can be described in terms of {\em closed families} of $B_K$-groups. \par
Moreover, each ideal $\mathsf{e}_{L,\varphi}$ associated to a $B_K$-group $(L,\varphi)$ has a unique maximal proper subideal $\mathsf{e}_{L,\varphi}^0$. The quotient $S_{L,\varphi}=\mathsf{e}_{L,\varphi}/\mathsf{e}_{L,\varphi}^0$ is a simple $\FB_K$-module. The evaluations of this simple module can be precisely described, as well as its minimal groups, and this yields a new example of a simple module over a Green biset functor with several isomorphism classes of minimal groups.\par
Finally, when $p$ is a prime number, we consider the restriction $\FB_K^{(p)}$ of $\FB_K$ to finite $p$-groups, and we describe completely the lattice of ideals of this Green biset functor. We show in particular that this lattice is always finite. As a byproduct, we get some examples of Green $p$-biset functors without non zero proper ideals.
\section{Review of shifted Green biset functors}
We quickly recall some definitions and basic notions on biset functors for finite groups, and refer to~\cite{bisetfunctorsMSC} for details. Let $\F$ be a field of characteristic 0. The biset category $\F\mathcal{C}$ of finite groups has all finite groups as objects. If $G$ and $H$ are finite groups, then $\Hom_{\F\mathcal{C}}(G,H)=\F\otimes_\Z B(H,G)$, where $B(H,G)$ is the Grothendieck group of finite $(H,G)$-bisets. Composition in $\F\mathcal{C}$ is induced by the product $(V,U)\mapsto V\times_HU=(V\times U)/H$, where $V$ is a $(K,H)$-biset and $U$ a $(H,G)$-biset, and $H$ acts on $(V\times U)$ by $(v,u)\cdot h=(vh,h^{-1}u)$. A biset functor over $\F$ is an $\F$-linear functor from $\F\mathcal{C}$ to the category of $\F$-vector spaces.\par
Any biset is a disjoint union of transitive ones, and any transitive $(H,G)$-biset is of the form $(H\times G)/L$, where $L$ is a subgroup of $(H\times G)$. Denoting by $p_1:H\times G\to H$ and $p_2:H\times G\to G$ the first and second projections, we set $k_1(L)=p_1(L\cap\Ker\,p_2)$ and $k_2(L)=p_2(L\cap\Ker\,p_1)$. The biset $(H\times G)/L$ factors as the composition
$$(H\times G)/L\cong \Ind_{p_1(L)}^H\circ\Inf_{p_1(L)/k_1(L)}^{p_1(L)}\circ\Iso(\alpha)\circ\Def_{p_2(L)/k_2(L)}^{p_2(L)}\circ\Res_{p_2(L)}^G$$
of elementary bisets called {\em induction}, {\em inflation}, {\em isomorphism}, {\em deflation}, and {\em restriction}, where $\alpha:p_2(L)/k_2(L)\to p_1(L)/k_1(L)$ is the canonical isomorphism sending $bk_2(L)$ to $ak_1(L)$ for $(a,b)\in L$. These elementary morphisms generate all morphisms in the category $\F\mathcal{C}$.\par
A Green biset functor $A$ over $\F$ (cf. Section~8.5 of~\cite{bisetfunctorsMSC}) is a biset functor with additional bilinear products $A(G)\times A(H)\to A(G\times H)$, denoted by $(\alpha,\beta)\mapsto \alpha\times \beta$, which are associative and bifunctorial. There is also an identity element $\varepsilon_A\in A(\un)$. \par
A left $A$-module $M$ is then defined similarly as a biset functor with products $A(G)\times M(H)\to M(G\times H)$ which are associative, bifunctorial, and unital. Left $A$-modules form an abelian category denoted by $\AMod{A}$. A left ideal of $A$ is an $A$-submodule of the left $A$-module $A$. \par
When $A$ is a Green functor, each evaluation $A(G)$ is an $\F$-algebra for the product
$$\alpha,\beta\in A(G)\mapsto \alpha\cdot\beta=A\big(\Iso(\delta)\circ\Res_\Delta^{G\times G}\big)(\alpha\times\beta)\mvirg$$
where $\Delta$ is the diagonal subgroup of $G\times G$, and $\delta:\Delta\to G$ the canonical isomorphism. The identity element of this algebra is $A(\Inf_\un^G)(\varepsilon_A)$. If $M$ is an $A$-module, each evaluation $M(G)$ is endowed with an $A(G)$-module structure defined similarly. By Proposition 2.16 of~\cite{tounkara-jofa}, a biset subfunctor $I$ of $A$ is an ideal if and only if $I(G)$ is an ideal of the algebra $A(G)$, for any finite group~$G$.\par
A Green biset functor $A$ is called {\em commutative} (cf. \cite{centros} for details) if the algebra $A(G)$ is commutative, for any $G$.\par
A fundamental example of Green biset functor is the Burnside functor sending a finite group $G$ to $\FB(G)=\FB(G,1)$, where $B(G)$ is the Burnside group of $G$. The products maps $\FB(G)\times\FB(H)\to\FB(G\times H)$ are induced by the cartesian product sending a $G$-set $X$ and an $H$-set $Y$ to the $(G\times H)$-set $X\times Y$. An $\FB$-module is precisely a biset functor over $\F$.\par
Let $K$ be a finite group. A Green biset functor $A$ over $\F$ can be {\em shifted} by $K$. This gives a new Green biset functor $A_K$ defined for a finite group $G$ by
$$A_K(G)=A(G\times K)\mpoint$$
For finite groups $G$ and $H$ and a finite $(H,G)$-biset $U$, the map 
$$A_K(U):A_K(G)\to A_K(H)$$
is the map $A(U\times K)$, where $U\times K$ is viewed as a $(H\times K,G\times K)$-biset in the obvious way, letting $K$ act on both sides on $U\times K$ by multiplication on the second component. For an arbitrary element $\alpha\in  \FB(H,G)$, that is an $\F$-linear combination of $(H,G)$-bisets, the map $ A_K(\alpha): A_K(G)\to A_K(H)$ is defined by $\F$-linearity.\par
This endows $ A_K$ with a biset functor structure. Moreover, for finite groups $G$ and $H$, the product
$$\times_{A_K}:  A_K(G)\times A_K(H)\to  A_K(G\times H)$$
is defined as follows: if $\alpha\in A_K(G)=A(G\times K)$ and $\beta\in A_K(H)=A(H\times K)$, then $\alpha\times\beta\in A(G\times K\times H\times K)$. We set
$$\alpha\times_{A_K}\beta=A\big(\Iso(\delta)\circ \Res_{\Delta}^{G\times K\times H\times K}\big)(\alpha\times\beta)\mvirg$$
where $\Delta=\{(g,k,h,k)\mid g\in G, h\in H, k\in K\}$, and $\delta$ is the isomorphism $\Delta\to G\times H\times K$ sending $(g,k,h,k)$ to $(g,h,k)$. The identity element $\varepsilon_{A_K}$ is $A(\Inf_\un^K)(\varepsilon_K)$.\par
For a finite group $G$, the algebra structure on $A_K(G)$ is simply the algebra structure on $A(G\times K)$ defined for the Green functor $A$.\par 
All these notion can be extended to functors from an admissible subcategory $\mathcal{D}$ of the biset category (cf. Chapter~4 of~\cite{bisetfunctorsMSC}), which is moreover closed under taking direct products of finite groups. We have then the notions of $\mathcal{D}$-biset functors and $\mathcal{D}$-Green biset functors, as well as modules over them.\par
In this paper, we will consider the shifted Burnside functor $\FB_K$, and its restriction $\FB_K^{(p)}$ to finite $p$-groups, for a prime $p$. 
A fundamental classical result is that for any finite group $G$, the algebra $\FB_K(G)$ is a split semisimple commutative algebra, with primitive idempotents $e_L^{G\times K}$ indexed by subgroups $L$ of $G\times K$, up to conjugation. The explicit formula for $e_L^{G\times K}$, due to Gluck (\cite{gluck}) and Yoshida~(\cite{yoshidaidemp}) is
$$e_L^{G\times K}=\frac{1}{|N_{G\times K}(L)|}\sum_{X\leq L}|X|\,\mu(X,L)\,[(G\times K)/X]\mvirg$$
where $X$ runs through all subgroups of $L$, where $\mu$ is the M\"obius function of the poset of subgroups of $G\times K$, and $[(G\times K)/X]$ is the isomorphism class of the transitive $(G\times K)$-set $(G\times K)/X$.\par
\begin{mth}{Notation} When $N$ is a normal subgroup of a finite group $L$, let
$$m_{L,N}=\frac{1}{|L|}\sum_{\substack{X\leq L\\XN=L}}|X|\mu(X,L)\mpoint$$
\end{mth}
\begin{mth}{Lemma} \label{deflation}Let $G$ be a finite group, and $L$ be a subgroup of $G\times K$. If $N$ is a normal subgroup of $G$, then 
$$\FB_K\big(\Def_{G/N}^G\big)(e_L^{G\times K})=\lambda \,m_{L,L\cap(N\times\un)}\,e_{\sur{L}}^{(G/N)\times K}\mvirg$$
where $\sur{L}$ is the image of $L$ by the projection $G\times K\to (G/N)\times K$, and $\lambda=\frac{|N_{(G/N)\times K}(\sur{L}):\sur{L}|}{|N_{G\times K}(L):L|}$.
\end{mth}
\pf Indeed 
$$\FB_K\big(\Def_{G/N}^G\big)(e_L^{G\times K})=\FB\big(\Def_{(G\times K)/(N\times \un)}^{G\times K}\big)(e_L^{G\times K})\mpoint$$
The result now follows from Assertion 4 of Theorem 5.2.4 of~\cite{bisetfunctorsMSC}.\endpf

\section{Ideals generated by idempotents}
We now introduce a category $\grpover{K}$, similar to the {\em comma category} over~$K$: its objects are the same, but morphisms are slightly different.
\begin{mth}{Definition} \label{grpover}\begin{itemize}
\item For a finite group $K$, let $\grpover{K}$ denote the following category:
\begin{itemize}
\item The objects are {\em finite groups over $K$}, i.e. pairs $(L,\varphi)$, where $L$ is a finite group and $\varphi:L\to K$ is a group homomorphism. 
\item A morphism $f:(L,\varphi)\to (L',\varphi')$ of groups over $K$ in the category $\grpover{K}$ is a group homomorphism $f:L\to L'$ such that there exists some inner automorphism $i$ of $K$ with $i\circ\varphi=\varphi'\circ f$.
\item The composition of morphisms in $\grpover{K}$ is the composition of group homomorphisms, and the identity morphism of $(L,\varphi)$ is the identity automorphism of $L$.
\end{itemize}
\item If $(L,\varphi)$ and $(L',\phi')$ are groups over $K$, we say that $(L',\varphi')$ is a {\em quotient} of $(L,\varphi)$, and we note $(L,\varphi)\twoheadrightarrow (L',\varphi')$, if there exists a morphism $f\in\Hom_{\grpover{K}}\big((L,\varphi),(L',\varphi')\big)$ with $f:L\to L'$ surjective. In this case, we will say that $f$ is a {\em surjective morphism} from $(L,\varphi)$ to $(L',\varphi')$.
\end{itemize}
\end{mth} 
\begin{rem}{Remarks} 
\begin{enumerate}
\item Using the well known fact that the epimorphisms in the category of (finite) groups are the surjective group homomorphisms (cf.~\cite{cwm} I.5 Exercise 5), one can show that a morphism $f\in\Hom_{\grpover{K}}\big((L,\varphi),(L',\varphi')\big)$ is an epimorphism in $\grpover{K}$ if and only if $f:L\to L'$ is surjective, that is, if $f$ is a surjective morphism. We will not use this fact here, except as a motivation to the use of the word ``quotient'' in Definition~\ref{grpover}.
\item A morphism $f:(L,\varphi)\to(L',\varphi')$ in $\grpover{K}$ is an isomorphism if and only if $f:L\to L'$ is an isomorphism of groups.
\item If $(L',\varphi')$ is a quotient of $(L,\varphi)$, and if $(L,\varphi)$ is a quotient of $(L',\varphi')$, then $(L,\varphi)$ and $(L',\varphi')$ are isomorphic in $\grpover{K}$. Indeed any surjective morphism from $(L,\varphi)$ to $(L',\varphi')$ is an isomorphism, for $L$ and $L'$ have the same order.
\item Clearly, the relation ``being quotient of'' on the class of groups over $K$ is transitive. In particular, any group over $K$ isomorphic in $\grpover{K}$ to a quotient of $(L,\varphi)$ is itself a quotient of $(L,\varphi)$, and also a quotient of any group over~$K$ isomorphic to $(L,\varphi)$ in $\grpover{K}$.
\end{enumerate}
\end{rem}
\begin{mth}{Notation} When $(L,\varphi)$ is a group over $K$, we denote by $L_\varphi$ the subgroup of $L\times K$ defined by
$$L_\varphi=\{\big(l,\varphi(l)\big)\mid l\in L\}\mpoint$$
\end{mth}
\begin{mth}{Theorem} \label{reduction1}Let $I$ be an ideal of the Green biset functor $\FB_K$.
 If $G$ is a finite group and $L$ is a subgroup of $G\times K$, the following conditions are equivalent:
\begin{enumerate}
\item The idempotent $e_{L}^{G\times K}$ belongs to $I(G)$.
\item The idempotent $e_{L_{p_2}}^{L\times K}$ belongs to $I(L)$, where $p_2:L\to K$ is the restriction to $L$ of the second projection homomorphism $G\times K\to K$.
\end{enumerate}
\end{mth}
\pf \fbox{$1\Rightarrow 2$} Let $\widehat{L}=L_{p_1}\subseteq L\times G$, where $p_1:L\to G$ is the restriction to $L$ of the first projection homomorphism $G\times K\to G$. Thus
$$p_1(\widehat{L})=L,\;k_1(\widehat{L})=\un\times k_2(L),\;p_2(\widehat{L})=p_1(L),\;k_2(\widehat{L})=\un\mpoint$$
It follows that the $(L,G)$-biset $U=(L\times G)/\widehat{L}$ factors as
$$U\cong \Inf_{L/N}^L\circ\Iso(\theta^{-1})\circ\Res_{p_1(L)}^G\mvirg$$
where $N=\un\times k_2(L)\normal L$ and $\theta:L/N\to p_1(G)$ is the canonical isomorphism induced by the first projection $p_1:L\to G$.\par
Now if $e_{L}^{G\times K}$ belongs to $I(G)$, its restriction $\FB_K(\Res_{p_1(L)}^G)(e_{L}^{G\times K})$ belongs to $I(G)$. But
\begin{align*}
\FB_K(\Res_{p_1(L)}^G)(e_{L}^{G\times K})&=\FB\big(\Res_{p_1(L)\times K}^{G\times K}\big)(e_{L}^{G\times K})\\
&=\sum_{L'}e_{L'}^{p_1(L)\times K},
\end{align*}
where $L'$ runs through a set of representatives of $\big(p_1(L)\times K\big)$-conjugacy classes of subgroups of $p_1(L)\times K$ which are conjugate to $L$ in $G\times K$ (cf. \cite{bisetfunctorsMSC}, Theorem 5.2.4, Assertion 1). In particular, the group $L$ is one of them, and 
$$e_L^{p_1(L)\times K}\cdot \FB_K(\Res_{p_1(L)}^G)(e_{L}^{G\times K})=e_L^{p_1(L)\times K}\in I\big(p_1(L)\big)\mpoint$$
It follows that $\FB_K\big(\Iso(\theta^{-1})\big)(e_L^{p_1(L)\times K})\in I(L/N)$. \par
But $\FB_K\big(\Iso(\theta^{-1})\big)=\FB\big(\Iso(\theta_K^{-1})\big)$, where $\theta_K=\theta\times \Id_K$ is the isomorphism from $(L/N)\times K$ to $p_1(L)\times K$ deduced from $\theta$. It follows that $e_{\sur{L}}^{(L/N)\times K}\in I(L/N)$, where $\sur{L}=\theta_K^{-1}(L)=\{\big(lN,p_2(l)\big)\mid l\in L\}$. Now 
\begin{align*}
\FB_K(\Inf_{L/N}^L)(e_{\sur{L}}^{(L/N)\times K})&=\FB\big(\Inf_{(L/N)\times K}^{L\times K}\big)(e_{\sur{L}}^{(L/N)\times K})\\
&=\sum_Xe_X^{L\times K}\in I(L)\mvirg
\end{align*}
where $X$ runs through a set of $(L\times K)$-conjugacy classes of subgroups of $L\times K$ which map to a conjugate of $\sur{L}$ through the surjection $L\times K\to (L/N)\times K$ (cf. \cite{bisetfunctorsMSC}, Theorem 5.2.4, Assertion 3).\par
The group $L_{p_2}$ is one of these subgroups, hence
$$e_{L_{p_2}}^{L\times K}\cdot\FB_K(\Inf_{L/N}^L)(e_{\sur{L}}^{(L/N)\times K})=e_{L_{p_2}}^{L\times K}\in I(L)\mvirg$$
as was to be shown.\mpn
\fbox{$2\Rightarrow 1$} We now consider the opposite $(G,L)$-biset $U\op\cong (G\times L)/\widetilde{L}$, where $\widetilde{L}=\{\big(p_1(l),l\big)\mid l\in L\}$, which factors as
$$U\op\cong \Ind_{p_1(L)}^G\circ\Iso(\theta)\circ\Def_{L/N}^L\mpoint$$
If $e_{L_{p_2}}^{L\times K}\in I(L)$, then $u=\FB_K(U\op)(e_{L_{p_2}}^{L\times K})$ belongs to $I(G)$. By Lemma~\ref{deflation}
$$\FB_K\big(\Def_{L/N}^L\big)(e_{L_{p_2}}^{L\times K})=\lambda\,m_{L_{l_2},L_{p_2}\cap(N\times\un)}\,e_{\sur{L_{p_2}}}^{(L/N)\times K}\mvirg$$
where $\sur{L_{p_2}}$ is the image of $L_{p_2}$ by the projection $L\times K\to (L/N)\times K$, and $\lambda$ is some non zero rational number.
Now the intersection
$$L_{p_2}\cap(N\times\un)=\{\big(a,b),b\big)\mid (a,b)\in L\}\cap \Big(\big(\un\times k_2(L)\big)\times\un\Big)$$
is trivial. It follows that $m_{L_{l_2},L_{p_2}\cap(N\times\un)}=1$, and
\begin{align*}
u&=\lambda\,\FB_K\big(\Ind_{p_1(L)}^G\circ\Iso(\theta)\big)(e_{\sur{L_{p_2}}}^{(G/N)\times K)})\\
&=\lambda\,\FB\big(\Ind_{p_1(L)\times K}^{G\times K}\circ\Iso(\theta_K)\big)(e_{\sur{L_{p_2}}}^{(G/N)\times K)})\mpoint
\end{align*}
Now for $(a,b)\in L$, the image by $\theta_K=\theta\times\Id_K$ of $\big((a,b),b\big)(N\times\un)\in\sur{L_{p_2}}$ is the element $\big(p_1(a,b),b\big)=(a,b)$ of $p_1(L)\times K$. Hence $\theta_K(\sur{L_{p_2}})$ identifies with~$L$, viewed as a subgroup of $p_1(L)\times K$, and
$$u=\lambda\,\FB\big(\Ind_{p_1(L)\times K}^{G\times K})(e_L^{p_1(L)\times K})=\lambda\lambda'e_L^{G\times K}\mvirg$$
for some non zero rational number $\lambda'$ (cf. \cite{bisetfunctorsMSC}, Theorem 5.2.4, Assertion 2). 
Since $u\in I(G)$ and $\lambda\lambda'\neq 0$, it follows that $e_{L}^{G\times K}\in I(G)$, as was to be shown.\endpf
\begin{mth}{Corollary} \label{reduction} Let $G$ be a finite group, and $L$ be a subgroup of $G\times\nolinebreak K$. Then the ideal of $\FB_K$ generated by $e_L^{G\times K}$ is equal to the ideal of $\FB_K$ generated by $e_{L_{p_2}}^{L\times K}$
\end{mth}
\pf Indeed, denoting by $I$ the ideal generated by $e_L^{G\times K}$, and by $J$ the ideal generated by $e_{L_{p_2}}^{L\times K}$, we have
\begin{align*}
e_L^{K\times G}\in I(G) &\Rightarrow e_{L_{p_2}}^{L\times K}\in I(L)\Rightarrow J\subseteq I\mvirg\\
e_{L_{p_2}}^{L\times K}\in J(L)&\Rightarrow e_L^{G\times K}\in J(G)\Rightarrow I\subseteq J\mvirg
\end{align*}
so $I=J$.\endpf
\begin{mth}{Notation} Let $(L,\varphi)$ be a group over $K$. We denote by $\mathsf{e}_{L,\varphi}$ the ideal of $\FB_K$ generated by $e_{L_\varphi}^{L\times K}\in\FB_K(L)$.
\end{mth}
\begin{mth}{Lemma} \label{quotient} Let $(L,\varphi)$ and $(M,\psi)$ be groups over $K$. 
\begin{enumerate}
\item If $(M,\psi)\twoheadrightarrow (L,\varphi)$, then $\mathsf{e}_{M,\psi}\subseteq \mathsf{e}_{L,\varphi}$.
\item In particular, if $(M,\psi)$ is isomorphic to $(L,\varphi)$, then $\mathsf{e}_{M,\psi}= \mathsf{e}_{L,\varphi}$.
\end{enumerate}
\end{mth}
\pf 1. Let $s:M\to L$ be a surjective group homomorphism, and $i$ be an inner automorphism of $K$ such that $i\circ\psi=\varphi\circ s$. Let $U$ denote the set~$L$, viewed as an $(M,L)$-biset for the action given by $m\cdot u\cdot l=s(m)ul$, for $m\in M$ and $u,l\in L$. There is an isomorphism of $(M,L)$-bisets
$$U\cong \Inf_{M/\Ker\,s}^M\circ\Iso(\alpha^{-1})\mvirg$$
where $\alpha:\sur{M}=M/\Ker\,s\to L$ is the group isomorphism induced by $s$.\par
Let $u=\FB_K(U)(e_{L_\varphi}^{L\times K})\in\mathsf{e}_{L,\varphi}(M)$. Then
$$u=\FB\big(\Inf_{\sur{M}\times K}^{M\times K}\circ\Iso(\alpha_K^{-1})\big)(e_{L_\varphi}^{L\times K})\mvirg$$
where $\alpha_K=\alpha\times \Id_K:\sur{M}\times K\to L\times K$. Then
$$\alpha_K^{-1}(L_\varphi)=\{\big(\alpha^{-1}(l),\varphi(l)\big)\mid l\in L\}=\{\big(m\Ker\,s,\varphi \circ s(m)\big)\mid m\in M\}\mpoint$$
It follows that $\FB\big(\Iso(\alpha_K^{-1})\big)(e_{L_\varphi}^{L\times K})=e_{\sur{M}_\theta}^{\sur{M}\times K}$, where $\theta:\sur{M}\to K$ is defined by $\theta(m\Ker\,s)=\varphi \circ s(m)$. In particular $e_{\sur{M}_\theta}^{\sur{M}\times K}\in\mathsf{e}_{L,\varphi}(\sur{M})$. Now
$$u=\FB\big(\Inf_{\sur{M}\times K}^{M\times K}\big)(e_{\sur{M}_\theta}^{\sur{M}\times K})=\sum_Xe_X^{M\times K}\mvirg$$
where $X$ runs through a set of representatives of conjugacy classes of subgroups of $M\times K$ such that the projection of $X$ in $\sur{M}\times K$ is conjugate to~$\sur{M}_\theta$. The subgroup $M_{\varphi\circ s}$ is one of these subgroups, so $e_{M_{\varphi\circ s}}^{M\times K}\cdot u$ is a non zero scalar multiple of $e_{M_{\varphi\circ s}}^{M\times K}$ lying in $\mathsf{e}_{L,\varphi}(M)$. Hence $e_{M_{\varphi\circ s}}^{M\times K}\in \mathsf{e}_{L,\varphi}(M)$. Now $\varphi\circ s=i\circ\psi$, where $i$ is an inner automorphism of $K$. This implies readily that the subgroups $M_{i\circ\psi}$ and $M_\psi$ of $M\times K$ are conjugate. It follows that 
$$e_{M_\psi}^{M\times K}=e_{M_{i\circ\psi}}^{M\times K}=e_{M_{\varphi\circ s}}^{M\times K}\in \mathsf{e}_{L,\varphi}(M)\mvirg$$
that is $\mathsf{e}_{M,\psi}\subseteq \mathsf{e}_{L,\varphi}$, proving Assertion 1.\mpn
Now if $f:(L,\varphi)\to(M,\psi)$ is an isomorphism in $\grpover{K}$, the group homomorphism $f:M\to L$ is an isomorphism. Then $(M,\psi)$ and $(L,\varphi)$ are quotient of one another, so $\mathsf{e}_{M,\psi}= \mathsf{e}_{L,\varphi}$, proving Assertion 2.\endpf
\begin{mth}{Notation} We fix a set $\mathcal{S}_K$ of isomorphism classes of objects in the category $\grpover{K}$.
\end{mth}
\begin{mth}{Proposition} \label{engendrent}Let $I$ be an ideal of $\FB_K$. Then $I$ is equal to the sum of the ideals $\mathsf{e}_{L,\varphi}$ it contains. More precisely, if 
$$\mathcal{A}_I=\{(L,\varphi)\in\mathcal{S}_K\mid \mathsf{e}_{L,\varphi}\subseteq I\}\mvirg$$
we have $I=\sum_{(L,\varphi)\in\mathcal{A}_I}\limits \mathsf{e}_{L,\varphi}$. It follows that the ideals of $\FB_K$ form a set.
\end{mth}
\pf Let $J=\sum_{\substack{(L,\varphi)\in\mathcal{S}_K\\\mathsf{e}_{L,\varphi}\subseteq I}}\limits\mathsf{e}_{L,\varphi}$. Then obviously $J\subseteq I$. Moreover, if $(M,\psi)$ is a group over $K$ such that $\mathsf{e}_{M,\psi}\subseteq I$, then $\mathsf{e}_{M,\psi}\subseteq J$: indeed, there is some $(L,\varphi)\in\mathcal{S}_K$ isomorphic to $(M,\psi)$, and $\mathsf{e}_{M,\psi}=\mathsf{e}_{L,\varphi}$ by Lemma~\ref{quotient}. Conversely, let $G$ be a finite group, and $u\in I(G)$. Then $u$ is a linear combination
$$u=\sum_{L}\lambda_Le_L^{G\times K}$$
with coefficients $\lambda_L$ in $\F$, of idempotents $e_L^{G\times K}$, where $L$ runs through a set $S$ of representatives of conjugacy classes of subgroups of $G\times K$. Then for any $L\in S$, we have $e_L^{G\times K}\cdot u=\lambda_Le_L^{G\times K}\in I(G)$, hence $e_L^{G\times K}\in I(G)$ if $\lambda_L\neq 0$. So in this case, the ideal of $\FB_K$ generated by $e_L^{G\times K}$ is contained in $I$. This ideal is equal to $\mathsf{e}_{L,p_2}$, by Corollary~\ref{reduction}, thus $\mathsf{e}_{L,p_2}\subseteq J$ by the above observation. Hence $e_L^{G\times K}\in \mathsf{e}_{L,p_2}(G)\subseteq J(G)$. It follows that
$$u=\sum_{\substack{L\in S\\\lambda_L\neq 0}}\lambda_Le_L^{G\times K}$$
also belongs to $J(G)$. Hence $I(G)\subseteq J(G)$, so $I(G)=J(G)$ since $J\subseteq I$. As~$G$ was arbitrary, it follows that $I=J$.\par
Now an ideal $I$ of $\FB_K$ is determined by the subset $\mathcal{A}_I$ of $\mathcal{S}_K$, so the class of ideals of $\FB_K$ is in one to one correspondence with a set of subsets of $\mathcal{S}_K$. Hence this class is a set.
\endpf
\begin{mth}{Lemma} \label{premier}Let $\mathcal{A}$ be a set of ideals of $\FB_K$, and $(M,\psi)$ be a group over~$K$. The following are equivalent:
\begin{enumerate}
\item $\mathsf{e}_{M,\psi}\subseteq \sum_{I\in\mathcal{A}}\limits I$.
\item There exists $I\in\mathcal{A}$ such that $\mathsf{e}_{M,\psi}\subseteq I$.
\end{enumerate}
\end{mth}
\pf Clearly 2 implies 1. Now 1 is equivalent to saying that 
$$e_{M_\psi}^{M\times K}\in \sum_{I\in\mathcal{A}}\limits I(M)\mpoint$$
If this holds, there exists $I\in\mathcal{A}$ and $u\in I(M)$ such that $e_{M_\psi}^{M\times K}\cdot u\neq 0$. Now $e_{M_\psi}^{M\times K}\cdot u\in I(M)$, and moreover there is a scalar $\lambda\in\F$ such that $e_{M_\psi}^{M\times K}\cdot u=\lambda e_{M_\psi}^{M\times K}\neq 0$. Hence $\lambda\neq 0$, and $e_{M_\psi}^{M\times K}\in I(M)$. In other words $\mathsf{e}_{M,\psi}\subseteq I$, so 1 implies 2.\endpf 
\section{$B_K$-groups}
In view of Proposition~\ref{engendrent}, every ideal of $\FB_K$ is a sum of ideals $\mathsf{e}_{L,\varphi}$, where $(L,\varphi)$ runs in a some subset of $\mathcal{S}_K$. In view of Lemma~\ref{premier}, to describe the inclusions between such sum of ideals $\mathsf{e}_{L,\varphi}$, it suffices to describe elementary inclusions of the form $\mathsf{e}_{M,\psi}\subseteq\mathsf{e}_{L,\varphi}$, where $(L,\varphi)$ and $(M,\psi)$ are groups over~$K$. Lemma~\ref{quotient} shows that it is the case if $(M,\psi)\twoheadrightarrow (L,\varphi)$. Moreover:
\begin{mth}{Theorem} \label{surjection forte}Let $s:(M,\psi)\to (L,\varphi)$ be a surjective morphism in~$\grpover{K}$. If $m_{M,\Ker\,s}\neq 0$, then $\mathsf{e}_{M,\psi}=\mathsf{e}_{L,\varphi}$.
\end{mth}
\pf We already know from Lemma~\ref{quotient} that $\mathsf{e}_{M,\psi}\subseteq\mathsf{e}_{L,\varphi}$, so it suffices to prove the reverse inclusion. We first observe that since there exists an inner automorphism $i$ of $K$ such that $i\circ\psi=\varphi\circ s$, we have $\Ker\,s\leq \Ker\,(i\circ\psi)=\Ker\,\psi$. So there is a group homomorphism $\sur{\psi}:\sur{M}=M/\Ker\,s\to K$ such that $\psi=\sur{\psi}\circ\pi$, where $\pi:M\to\sur{M}$ is the projection map.\par
Now let $V$ be the set $L$, viewed as an $(L,M)$-biset for the action defined by $l\cdot v\cdot m=lvs(m)$, for $l,v\in V$ and $m\in M$ (in other words $V=U\op$, where $U$ is the $(M,L)$-biset introduced in the proof of Lemma~\ref{quotient}). Then there is an isomorphism of $(L,M)$-bisets
$$V\cong \Iso(\alpha)\circ\Def_{M/\Ker\,s}^M\mvirg$$
where $\alpha:\sur{M}\to L$ is the group isomorphism induced by $s$, i.e. such that $s=\alpha\circ\pi$.\par
Let $v=\FB_K(V)(e_{M_\psi}^{M\times K})\in \mathsf{e}_{M,\psi}(L)$. By Lemma~\ref{deflation}
$$\FB_K\big(\Def_{M/\Ker\,s}^M\big)(e_{M_\psi}^{M\times K})=\lambda\,m_{M_\psi,M_\psi\cap(\Ker\,s\times\un)}\,e_{\sur{M_\psi}}^{\sur{M}\times K}\mvirg$$
where $\sur{M_\psi}$ is the image of $M_\psi$ by the projection $M\times K\to \sur{M}\times K$, and $\lambda$ is some non zero rational number. Then
$$v=\lambda\,\FB_K\big(\Iso(\alpha)\big)(e_{\sur{M_\psi}}^{\sur{M}\times K})=\lambda\,\FB\big(\Iso(\alpha_K))(e_{\sur{M_\psi}}^{\sur{M}\times K})\mvirg$$
where $\alpha_K=\alpha\times \Id_K:\sur{M}\times K\to L\times K$.
The image of $\sur{M}_{\sur{\psi}}$ under $\alpha_K$ is the subgroup 
$$\alpha_K(\sur{M}_{\sur{\psi}})=\big\{\big(\alpha(\sur{m}),\sur{\psi}(\sur{m})\big)\mid\sur{m}\in\sur{M}\big\}=\big\{\big(l,\sur{\psi}\circ\alpha^{-1}(l)\big)\mid l\in L\big\}\mpoint$$
Moreover, we have a diagram
$$\xymatrix@C=3ex@R=3ex{
M\ar@/_/[dddr]_-\psi\ar[dr]^-{\pi}\ar@/^/[drrr]^-s&&\\
&\sur{M}\ar[rr]_-\alpha\ar[dd]^-{\sur{\psi}}&&L\ar[dd]_-\varphi\\
&&&\\
&K\ar[rr]^-i&&K
}
$$
where the two triangles and the outer ``square'' commute. It follows that 
$$\varphi\circ\alpha\circ\pi=\varphi\circ s=i\circ\psi=i\circ\sur{\psi}\circ\pi\mvirg$$
hence $\varphi\circ\alpha=i\circ\sur{\psi}$ since $\pi$ is surjective. Hence
$\sur{\psi}\circ\alpha^{-1}=i^{-1}\circ\varphi$, and $\alpha_K(\sur{M}_{\sur{\psi}})=L_{i^{-1}\circ\varphi}$.\par
It follows that $v=\lambda\,e_{L_{i^{-1}\circ\varphi}}^{L\times K}$, and moreover $e_{L_{i^{-1}\circ\varphi}}^{L\times K}=e_{L_\varphi}^{L\times K}$ since $L_{i^{-1}\circ\varphi}$ and $L_\varphi$ are conjugate in $L\times K$. Finally $v=\lambda\,e_{L_\varphi}^{L\times K}$, so $e_{L_\varphi}^{L\times K}\in \mathsf{e}_{M,\psi}(L)$, since $v\in \mathsf{e}_{M,\psi}(L)$ and $\lambda\neq 0$. In other words $\mathsf{e}_{L,\varphi}\subseteq \mathsf{e}_{M,\psi}$, and finally $\mathsf{e}_{L,\varphi}= \mathsf{e}_{M,\psi}$, as was to be shown.\endpf
\begin{mth}{Notation} When $(M,\psi)$ is a group over $K$, and $Q$ is a normal subgroup of $M$ with $Q\leq\Ker\,\psi$, let $\psi/Q:M/Q\to K$ be the group homomorphism defined by $\psi=(\psi/Q)\circ \pi$, where $\pi$ is the projection $M\to M/Q$.
\end{mth}
Thus for any group $(M,\psi)$ over $K$, if $Q$ is a normal subgroup of $M$ contained in $\Ker\,\psi$, we get a surjective morphism $\pi:(M,\psi)\to (M/Q,\psi/Q)$ in $\grpover{K}$, with $\Ker\,\pi=Q$. If moreover $m_{M,Q}\neq 0$, we have $\mathsf{e}_{M,\psi}=\mathsf{e}_{M/Q,\psi/Q}$. This motivates the following:
\begin{mth}{Definition} Let $(L,\varphi)$ be a group over $K$. We say that $(L,\varphi)$ is a {\em $B_K$-group}, or a {\em $B$-group relative to $K$}, if $m_{L,N}=0$ for every non-trivial normal subgroup $N$ of $L$ contained in $\Ker\,\varphi$.
\end{mth}
\begin{rem}{Examples} \label{exBK}\begin{enumerate}
\item If $\varphi:L\to K$ is injective, then $(L,\varphi)$ is a $B_K$-group.
\item On the other hand, if $K=\un$, then a group over $K$ is a pair $(L,\varphi)$, where $L$ is a finite group and $\varphi:L\to \un$ is the unique morphism. Moreover the category $\grpover{\un}$ clearly identifies with the usual category of finite groups. With this identification, a $B_\un$-group is just a $B$-group (cf. Section~7.2 of \cite{doublact}, or Chapter~5 of~\cite{bisetfunctorsMSC}).
\end{enumerate}
\end{rem}
\begin{mth}{Lemma} \label{BK-group iso}Let $(L,\varphi)$ be a $B_K$-group. If $(M,\psi)$ is a group over $K$, and $(M,\psi)$ is isomorphic to $(L,\varphi)$ in $\grpover{K}$, then $(M,\psi)$ is a $B_K$-group.
\end{mth}
\pf Since $(M,\psi)$ is isomorphic to $(L,\varphi)$ in $\grpover{K}$, there exists a group isomorphism $f:L\to M$ and an inner automorphism $i$ of $K$ such that $\psi\circ f=i\circ\varphi$. If $P$ is a normal subgroup of $M$ contained in $\Ker\,\psi$, then $f^{-1}(P)$ is a normal subgroup of $L$ contained in $\Ker\,\varphi$, because 
$$i\circ\varphi\big(f^{-1}(P)\big)=\psi\circ f\big(f^{-1}(P)\big)=\psi(P)=\un$$
 and $i$ is an automorphism. Moreover $m_{L,f^{-1}(P)}=m_{M,P}$. If $P$ is non trivial, then $f^{-1}(P)$ is non trivial, so $m_{L,f^{-1}(P)}=m_{M,P}=0$, as was to be shown.\endpf
\begin{mth}{Theorem} \label{BK-groupe}Let $(L,\varphi)$ be a group over $K$.
\begin{enumerate}
\item If $Q$ is a normal subgroup of $L$, contained in $\Ker\,\varphi$, and maximal such that $m_{L,Q}\neq 0$, then $(L/Q,\varphi/Q)$ is a $B_K$-group, quotient of $(L,\varphi)$.
\item If $(P,\psi)$ is a $B_K$-group which is quotient of $(L,\varphi)$, and if $N$ is a normal subgroup of $L$ contained in $\Ker\,\varphi$ and such that $m_{L,N}\neq 0$, then $(P,\psi)$ is a quotient of $(L/N,\varphi/N)$.
\item In particular, if $P$ and $Q$ are normal subgroups of $L$, contained in $\Ker\,\varphi$, and maximal such that $m_{L,P}\neq 0\neq m_{L,Q}$, then $(L/P,\varphi/P)$ and $(L/Q,\varphi/Q)$ are isomorphic in $\grpover{K}$.
\end{enumerate}
\end{mth}
\pf 1. Let $P/Q$ be a normal subgroup of $L/Q$ contained in $\Ker\,(\varphi/Q)=\Ker\,\varphi/Q$. Then
 $P$ is a normal subgroup of $L$, and $Q\leq P\leq\Ker\,\varphi$. If $P/Q\neq\un$, i.e. if $Q<P$, then by maximality of $Q$ and Proposition~5.3.1 of~\cite{bisetfunctorsMSC}
$$m_{L,P}=0=m_{L,Q}m_{L/Q,P/Q}\mpoint$$
Since $m_{L,Q}\neq 0$, it follows that $m_{L/Q,P/Q}=0$, so $(L/Q,\varphi/Q)$ is a $B_K$-group, quotient of $(L,\varphi)$.\mpn
2. Since $(P,\psi)$ is a quotient of $(L,\varphi)$, there exists a surjective group homomorphism $s:L\to P$ and an inner automorphism $i$ of $K$ such that $\psi\circ s=i\circ\varphi$. It follows that $M=\Ker\,s$ is a normal subgroup of $L$ contained in $\Ker\,(i\circ\varphi)=\Ker\,\varphi$.\par
We have a diagram
$$\xymatrix@C=3ex@R=3ex{
L\ar@/_/[dddr]_-\varphi\ar[dr]^-{\pi_M}\ar@/^/[drrr]^-s&&\\
&L/M\ar[rr]_-{\sur{s}}\ar[dd]^-{\varphi/M}&&P\ar[dd]_-\psi\\
&&&\\
&K\ar[rr]^-i&&K
}
$$
where the two triangles and the outer ``square'' commute, and $\sur{s}$ is an isomorphism, the map $\pi_M:L\to L/M$ being the projection. As in the proof of Theorem~\ref{surjection forte}, we have
$$\psi\circ\sur{s}\circ\pi_M=\psi\circ s=i\circ \varphi=i\circ(\varphi/M)\circ\pi_M\mvirg$$
so $\psi\circ\sur{s}=i\circ(\varphi/M)$ since $\pi_M$ is surjective. It follows that $\sur{s}$ is an isomorphism from $(L/M,\varphi/M)$ to $(P,\psi)$ in $\grpover{K}$, so $(L/M,\varphi/M)$ is a $B_K$-group by Lemma~\ref{BK-group iso}.\par
Now by Proposition~5.3.3 of~\cite{bisetfunctorsMSC}
$$m_{L,N}=\frac{1}{|L|}\sum_{YN=YM=L}|Y|\mu(Y,L)m_{L/M,(Y\cap N)M/M}\mpoint$$
In particular, if $m_{L,N}\neq 0$, there exists $Y\leq L$ such that $YN=YM=L$ and $m_{L/M,(Y\cap N)M/M}\neq 0$. But since $N\subseteq \Ker\,\varphi$, the group $(Y\cap N)M/M$ is a normal subgroup of $L/M$ contained in $\Ker\,(\varphi/M)=\Ker\,\varphi/M$. Then since $m_{L/M,(Y\cap N)M/M}\neq 0$ and since $(L/M,\varphi/M)$ is a $B_K$-group, we have $(Y\cap N)M/M=\un$, i.e. $Y\cap N\subseteq Y\cap M$.\par
Consider now the following diagram:
$$
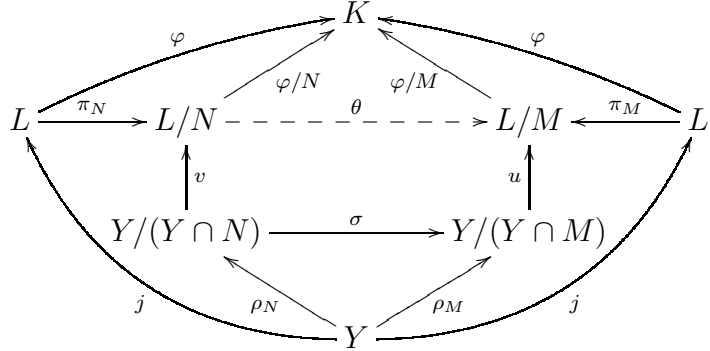
\xymatrix{
&&K&&\\
L\ar@/^/[urr]^-\varphi\ar[r]^-{\pi_N}&L/N\ar[ur]_-{\varphi/N}\ar@{-->}[rr]^-\theta&&L/M\ar[ul]^-{\varphi/M}&L\ar[l]_-{\pi_M}\ar@/_/[ull]_-\varphi\\
&Y/(Y\cap N)\ar[u]_-v\ar[rr]^-\sigma&&Y/(Y\cap M)\ar[u]^-u&\\
&&Y\ar@/^5ex/[lluu]^-j\ar[ul]^-{\rho_N}\ar[ur]_-{\rho_M}\ar@/_5ex/[uurr]_-j&&
}
$$
where 
\begin{itemize}
\item $j:Y\to L$ is the inclusion map, 
\item $\rho_N:Y\to Y/(Y\cap N)$ and $\rho_M:Y\to Y/(Y\cap M)$ are the projection maps,
\item $u:Y/(Y\cap M)\to L/M$ and $v:Y/(Y\cap N)\to L/N$ are the canonical isomorphisms $Y/(Y\cap M)\cong YM/M=L/M$ and $Y/(Y\cap N)\cong YN/N=L/N$, respectively,
\item $\sigma:Y/(Y\cap N)\to Y/(Y\cap M)$ is the projection map (as $Y\cap N\subseteq Y\cap M$),
\item $\theta:L/N\to L/M$ is defined as $\theta=u\circ\sigma\circ v^{-1}$. In particular $\theta$ is surjective.
\end{itemize}
We have $\pi_N\circ j=v\circ \rho_N$, since for any $y\in Y$
$$\pi_N\circ j(y)=\pi_N(y)=yN=v\big(y(Y\cap N)\big)=v\circ\rho_N(y)\mpoint$$
Similarly $\pi_M\circ j=u\circ\rho_M$. We also have $\sigma\circ \rho_N=\rho_M$. Then
$$\theta\circ\pi_N\circ j=\theta\circ v\circ\rho_N=u\circ\sigma\circ\rho_N=u\circ\rho_M=\pi_M\circ j\mpoint$$
Hence
$$(\varphi/M)\circ \theta\circ\pi_N\circ j=(\varphi/M)\circ \pi_M\circ j=\varphi\circ j=(\varphi/N)\circ\pi_N\circ j\mpoint$$
Since $\pi_N\circ j=v\circ\rho_N:Y\to L/N$ is surjective, it follows that
$$(\varphi/M)\circ \theta=(\varphi/N)\mpoint$$
Hence $\theta$ is a surjective morphism from $(L/N,\varphi/N)$ to $(L/M,\varphi/M)$ in $\grpover{K}$. As the latter is isomorphic to $(P,\psi)$ in $\grpover{K}$, it follows that $(P,\psi)$ is a quotient of $(L/N,\varphi/N)$, as was to be shown.\mpn
3. If $P$ and $Q$ are normal subgroups of $L$, contained in $\Ker\,\varphi$, and maximal such that $m_{L,P}\neq 0\neq m_{L,Q}$, then $(L/P,\varphi/P)$ and $(L/Q,\varphi/Q)$ are both $B_K$-groups by Assertion 1, and they are quotient of one another by Assertion 2. Hence they are isomorphic in $\grpover{K}$.\endpf
\begin{mth}{Notation} Let $(L,\varphi)$ be a group over $K$. If $Q$ is a normal subgroup of $L$, contained in $\Ker\,\varphi$, and maximal such that $m_{L,Q}\neq 0$, we denote by $\beta_K(L,\varphi)$ the quotient $(L/Q,\varphi/Q)$ of $(L,\varphi)$.
\end{mth}
\begin{rem}{Remark} As observed in Example~\ref{exBK}, when $K$ is trivial, a $B_K$-group is simply a $B$-group. Moreover, for any finite group $L$, if $u:L\to\un$ is the unique group homomorphism, then $\beta_\un(L,u)=\beta(L)$.
\end{rem}
The following corollary shows that $\beta_K(L,\varphi)$ is the largest $B_K$-group quotient of $(L,\varphi)$:
\begin{mth}{Corollary} \label{betaK}Let $(L,\varphi)$ be a group over $K$. 
\begin{enumerate}
\item $\beta_K(L,\varphi)$ is well defined up to isomorphism in $\grpover{K}$.
\item $\beta_K(L,\varphi)$ is a $B_K$-group, quotient of $(L,\varphi)$.
\item If $(P,\psi)$ is a $B_K$-group, quotient of $(L,\varphi)$, then $(P,\psi)$ is a quotient of $\beta_K(L,\varphi)$.
\item $\mathsf{e}_{L,\varphi}=\mathsf{e}_{\beta_K(L,\varphi)}$.
\end{enumerate}
\end{mth}
\pf 1. This follows from Assertion 3 of Theorem~\ref{BK-groupe}.\mpn
2. This follows from Assertion 1 of Theorem~\ref{BK-groupe}.\mpn
3. This follows from Assertion 2 of Theorem~\ref{BK-groupe}.\mpn
4. This follows from Theorem~\ref{surjection forte}, by definition of $\beta_K(L,\varphi)$
\endpf
\begin{mth}{Corollary} \label{betaK mGN}Let $s:(M,\psi)\twoheadrightarrow (L,\varphi)$ be a surjective morphism in $\grpover{K}$. Then $\beta_K(M,\psi)\cong \beta_K(L,\varphi)$ if and only if $m_{M,\Ker\,s}\neq 0$.
\end{mth}
\pf Indeed $\beta_K(L,\varphi)$ is a quotient of $(M,\psi)$, as it is a quotient of $(L,\varphi)$ and $s$ is surjective. Hence $\beta_K(L,\varphi)$ is a quotient of $\beta_K(M,\psi)$. Set $N=\Ker\,s$, so that $(L,\varphi)\cong(M/N,\psi/N)$.\par
 If $m_{M,N}\neq 0$, then since $\beta_K(M,\psi)$ is a $B_K$-group quotient of $(M,\psi)$, Assertion~2 of Theorem~\ref{BK-groupe} implies that $\beta_K(M,\psi)$ is a quotient of $(M/N,\psi/N)\cong(L,\varphi)$, hence of $\beta_K(L,\varphi)$. It follows that $\beta_K(M,\psi)\cong \beta_K(L,\varphi)$, as they are quotient of one another.\par
Conversely, suppose that $\beta_K(M,\psi)\cong \beta_K(L,\varphi)$, and let $P/N$ be a normal subgroup of $M/N$ contained in $\Ker\,(\psi/N)=\Ker\,\psi/N$ and maximal such that $m_{M/N,P/N}\neq 0$. Then the quotient $\big((M/N)\big/(P/N),(\psi/N)\big/(P/N)\big)\cong (M/P,\psi/P)$ is isomorphic to $\beta_K(M/N,P/N)\cong\beta_K(L,\varphi)$, hence to $\beta_K(M,\psi)$. Now if $Q$ is a normal subgroup of $M$ contained in $\Ker\,\psi$ and maximal such that $m_{M,Q}\neq 0$, the quotient $(M/Q,\psi/Q)$ is isomorphic to $\beta_K(M,\psi)\cong (M/P,\psi/P)$. In particular $M/Q\cong M/P$, and then $m_{M,P}=m_{M,Q}$ by Proposition 5.3.4 of~\cite{bisetfunctorsMSC}, so $m_{M,P}\neq 0$. But $m_{M,P}=m_{M,N}m_{M/N,P/N}$, so $m_{M,N}\neq 0$, as was to be shown.\endpf

\section{The ideals of $\FB_K$}
\begin{mth}{Notation and Definition} \begin{enumerate}
\item We let $\BgrK{K}$ denote the subset of $\mathcal{S}_K$ consisting of $B_K$-groups.
\item A subset $\mathcal{P}$ of $\BgrK{K}$ is said to be {\em closed} if
$$\forall (L,\varphi)\in\CP,\;\forall (M,\psi)\in \BgrK{K},\; (M,\psi)\twoheadrightarrow (L,\varphi)\;\implies \;(M,\psi)\in\CP\mpoint\vspace{-2ex}$$
\end{enumerate}
\end{mth}
\pagebreak[3]
\begin{mth}{Proposition} \label{PI}Let $I$ be an ideal of $\FB_K$, and
$$\CP_I=\{(L,\varphi)\in\BgrK{K}\mid \mathsf{e}_{L,\varphi}\subseteq I\}\mpoint$$
Then $\CP_I$ is a closed subset of $\BgrK{K}$, and $I=\sum_{(L,\varphi)\in\CP_I}\limits\mathsf{e}_{L,\varphi}$
\end{mth}
\pf The subset $\CP_I$ of $\BgrK{K}$ is closed by Lemma~\ref{quotient}. The second assertion follows from Proposition~\ref{engendrent} and Assertion 4 of Corollary~\ref{betaK}.\endpf
\begin{mth}{Theorem} \label{eLphi}Let $(L,\varphi)$ be a $B_K$-group. Then for any finite group~$G$
$$\mathsf{e}_{L,\varphi}(G)=\sum_X\F e_X^{G\times K}\mvirg$$
where $X$ runs through all subgroups of $G\times K$ such that $(X,p_2)\twoheadrightarrow (L,\varphi)$.
\end{mth}
\pf If $X\leq G\times K$ and $(X,p_2)\twoheadrightarrow (L,\varphi)$, then $\mathsf{e}_{X,p_2}\subseteq \mathsf{e}_{L,\varphi}$ by Lemma~\ref{quotient}. Equivalently $e_{X_{p_2}}^{X\times K}\in \mathsf{e}_{L,\varphi}(X)$, which is equivalent to $e_X^{G\times K}\in \mathsf{e}_{L,\varphi}(G)$, by Theorem~\ref{reduction1}. This proves that for each finite group $G$, the sum $E(G)=\sum_X\limits \F e_X^{G\times K}$, where $X\leq G\times K$ and $(X,p_2)\twoheadrightarrow (L,\varphi)$, is a subset of $\mathsf{e}_{L,\varphi}(G)$.\par
Moreover the map $\big(l,\varphi(l)\big)\in L_\varphi\mapsto l\in L$ is clearly an isomorphism $(L_\varphi,p_2)\to (L,\varphi)$ in $\grpover{K}$. In particular $(L_\varphi,p_2)\twoheadrightarrow (L,\varphi)$, and then by definition $e_{L_\varphi}^{L\times K}\in E(L)$. If we can prove that $G\mapsto E(G)$ defines an ideal $E$ of $\FB_K$, then we are done, because $E\subseteq \mathsf{e}_{L,\varphi}$ since $E(G)\subseteq \mathsf{e}_{L,\varphi}(G)$ for any~$G$, and $\mathsf{e}_{L,\varphi}\subseteq E$ because the generator $e_{L_\varphi}^{L\times K}$ of $\mathsf{e}_{L,\varphi}$ belongs to $E(L)$.\par
Since $E(G)$ is obviously an ideal of the algebra $\FB_K(G)$, for any $G$, all we have to do is to show that $E$ is a biset subfunctor of $\FB_K$, in other words that it is preserved by the elementary biset operations of induction, restriction, inflation, deflation, and transport by group isomorphism. For this, in what follows, we refer to Theorem~5.2.4 of~\cite{bisetfunctorsMSC}.\par
Let $X\leq G\times K$ be such that $(X,p_2)\twoheadrightarrow (L,\varphi)$, and suppose first that~$G$ is a subgroup of a group $H$. Then 
$$\FB_K\big(\Ind_G^H\big)(e_X^{G\times K})=\FB\big(\Ind_{G\times K}^{H\times K}\big)(e_X^{G\times K})=\lambda\,e_{X'}^{H\times K}$$
for some scalar $\lambda$, where $X'$ is the group $X$, viewed as a subgroup of $H\times K$. Clearly $(X',p_2)=(X,p_2)$, so $(X',p_2)\twoheadrightarrow (L,\varphi)$ and $e_{X'}^{H\times K}\in E(H)$. Hence $E$ is preserved by induction.\par
Assume now that $H$ is a subgroup of $G$. Then
$$\FB_K\big(\Res_H^G\big)(e_X^{G\times K})=\FB\big(\Res_{H\times K}^{G\times K}\big)(e_X^{G\times K})=\sum_Ye_Y^{H\times K}\mvirg$$
where $Y$ runs through a set of representatives of $(H\times K)$-conjugacy classes of subgroups of $H\times K$ which are conjugate to $X$ in $G\times K$. If $Y$ is such a subgroup, there exists $(g,k)\in G\times K$ such that $Y=X^{(g,k)}$. Then we have a commutative diagram
$$\xymatrix{
Y\ar[d]_-{p_2}\ar[r]^-\alpha&X\ar[d]^-{p_2}\\
K\ar[r]^-\beta&K
}
$$
where $\alpha$ is (left-)conjugation by $(g,k)$ and $\beta$ is (left-)conjugation by $k$. Since $\beta$ is an inner automorphism of $K$, and since $\alpha$ is a group isomorphism, it follows that $\alpha:(Y,p_2)\to (X,p_2)$ is an isomorphism in $\grpover{K}$. Hence $(Y,p_2)\twoheadrightarrow (L,\varphi)$, and $e_Y^{H\times K}\in E(H)$. It follows that $E$ is preserved by restriction.\par
Assume next that $G$ is a quotient of a group $H$ by a normal subgroup $N$. Then
$$\FB_K\big(\Inf_G^H\big)(e_X^{G\times K})=\FB\big(\Inf_{G\times K}^{H\times K}\big)(e_X^{G\times K})=\sum_Ye_Y^{H\times K}\mvirg$$
where $Y$ runs through a set of $(H\times K)$ conjugacy classes of subgroup of $H\times K$ which map to a conjugate of $X$ under the projection $\pi\times \Id_K:H\times K\to G\times K$, where $\pi:H\to G$ is the projection. Replacing $Y$ by a conjugate, which does not change $e_Y^{H\times K}$, we can assume that $Y$ is mapped to $X$ by $\pi\times \Id_K$. This gives a commutative diagram
$$\xymatrix@C=2ex{
Y\ar@{->>}[rr]^-{\pi\times\Id_K}\ar[dr]_-{p_2}&&X\ar[ld]^-{p_2}\\
&K&
}
$$
showing that $(Y,p_2)\twoheadrightarrow(X,p_2)$. Hence $(Y,p_2)\twoheadrightarrow (L,\varphi)$, so $e_Y^{H\times K}\in E(H)$, and $E$ is preserved by inflation.\par
As for deflation, we assume now that $H=G/N$, where $N\normal G$. Let $\pi:G\to H$ be the projection map. Then by Lemma~\ref{deflation}
$$\FB_K\big(\Def_H^G\big)(e_X^{G\times K})=\lambda\,m_{X,X\cap(N\times\un)}\,e_{\sur{X}}^{H\times K}\mvirg$$
where $\sur{X}$ is the image of $X$ under the projection $\pi\times\Id_K:G\times K\to H\times K$, and $\lambda$ is some non zero scalar. As above, we get a commutative diagram
$$\xymatrix@C=2ex{
X\ar@{->>}[rr]^-{s}\ar[dr]_-{p_2}&&\sur{X}\ar[ld]^-{p_2}\\
&K&
}
$$
where $s$ is the restriction of $\pi\times \Id_K$ to $X$. Then $s:(X,p_2)\to (\sur{X},p_2)$ is a surjective morphism in $\grpover{K}$. Setting $P=\Ker\,s=X\cap(N\times\un)$, we get an isomorphism $(\sur{X},p_2)\cong(X/P,p_2/P)$ in $\grpover{K}$.  Moreover $(L,\varphi)$ is a $B_K$-group quotient of $(X,p_2)$ by assumption. Then there are two cases: either $m_{X,P}=0$, and then $\FB_K\big(\Def_H^G\big)(e_X^{G\times K})=0\in E(H)$. Or $m_{X,P}\neq 0$, and then $(L,\varphi)$ is a quotient of $(X/P,p_2/P)\cong (\sur{X},p_2)$, by Assertion~2 of Theorem~\ref{BK-groupe}. It follows that $e_{\sur{X}}^{H\times K}\in E(H)$, so $\FB_K\big(\Def_H^G\big)(e_X^{G\times K})\in E(H)$ as well. This shows that $E$ is preserved by deflation.\par
Finally, it is clear that $E$ is preserved by group isomorphisms. This completes the proof of Theorem~\ref{eLphi}.\endpf
\begin{rem}{Remark} \label{base eLphi}Theorem~\ref{eLphi} implies that the set of idempotents $e_X^{G\times K}$, where $X$ runs through a set of {\em representatives of conjugacy classes} of subgroups of $G\times K$ such that $(X,p_2)\twoheadrightarrow (L,\varphi)$, is an {\em $\F$-basis} of $\mathsf{e}_{L,\varphi}(G)$.\vspace{-2ex}
\end{rem}
\begin{mth}{Corollary} \label{inclusion BK}Let $(L,\varphi)$ be a $B_K$-group, and $(M,\psi)$ be a group over~$K$. Then $\mathsf{e}_{M,\psi}\subseteq \mathsf{e}_{L,\varphi}$ if and only if $(M,\psi)\twoheadrightarrow (L,\varphi)$.
\end{mth}
\pf Indeed $\mathsf{e}_{M,\psi}\subseteq \mathsf{e}_{L,\varphi}$ if and only if $e_{M_\psi}^{M\times K}\in \mathsf{e}_{L,\varphi}(M)$, i.e. if and only if $(M_\psi,p_2)\twoheadrightarrow (L,\varphi)$. But we have already noticed at the beginning of the proof of Theorem~\ref{eLphi} that the map $\big(m,\psi(m)\big)\in M_\psi\mapsto m\in M$ is an isomorphism from $(M_\psi,p_2)$ to $(M,\psi)$ in $\grpover{K}$.\endpf
\begin{rem}{Remark} \label{cite centros}It was shown in Section 5.2.2 of~\cite{centros} that the category $\AMod{\FB_K}$ splits as a product
$$\AMod{\FB_K}\cong\prod_{H}=\AMod{e_H^K\FB_K}\mvirg\vspace{-2ex}$$
of categories of modules over smaller Green biset functors $e_H^K\FB_K$, where $H$ runs through a set of representatives of conjugacy classes of subgroups of~$K$. The functor $e_H^K\FB_K$ is the direct summand of $\FB_K$ obtained from the idempotent $e_H^K$ of $\FB_K(\un)\cong\FB(K)$. Its value at a group $G$ is the set of $\F$-linear combinations of idempotents $e_L^{G\times K}$ associated to subgroups $L$ for which $p_2(L)$ is conjugate to $H$ in $K$. This condition is equivalent to the existence of a surjective morphism $(L,p_2)\twoheadrightarrow (H,j_H)$, where $j_H:H\hookrightarrow K$ is the inclusion morphism. Since $(H,j_H)$ is a $B_K$-group by Example~\ref{exBK}, it follows that $e_H^K\FB_K=\mathsf{e}_{H,j_H}$.
\end{rem}
\begin{mth}{Theorem} \label{treillis}Let $\mathcal{I}_{\FB_K}$ be the lattice  of ideals of $\FB_K$, ordered by inclusion of ideals, and $\mathcal{C}l_{\BgrK{K}}$ be the lattice of closed subsets of $\BgrK{K}$, ordered by inclusion of subsets. Then the map
$$I\in\mathcal{I}_{\FB_K}\mapsto\CP_I=\{(L,\varphi)\in\BgrK{K}\mid \mathsf{e}_{L,\varphi}\subseteq I\}$$
is an isomorphism of lattices from $\mathcal{I}_{\FB_K}$ to $\mathcal{C}l_{\BgrK{K}}$. The inverse isomorphism is the map
$$\CP\in \mathcal{C}l_{\BgrK{K}}\mapsto I_\CP=\sum_{(L,\varphi)\in \CP}\mathsf{e}_{L,\varphi}\mpoint$$
In particular $\mathcal{I}_{\FB_K}$ is completely distributive.
\end{mth}
\pf  By Proposition~\ref{PI}, if $I$ is an ideal of $\FB_K$, then $\CP_I$ is a closed subset of $\BgrK{K}$, so the map $\alpha:I\mapsto \CP_I$ from $\mathcal{I}_{\FB_K}$ to $\mathcal{C}l_{\BgrK{K}}$ is well defined. It is moreover clearly order preserving. The map $\beta:\CP\mapsto \CP_I$ from $\mathcal{C}l_{\BgrK{K}}$ is also well defined and order preserving. By Proposition~\ref{PI} again, the composition $\beta\circ\alpha$ is the identity map of $\mathcal{I}_{\FB_K}$. Conversely, if $\CP\in \mathcal{C}l_{\BgrK{K}}$, then
$$\alpha\circ\beta=\big\{(M,\psi)\in\BgrK{K}\mid \mathsf{e}_{M,\psi}\subseteq \sum_{(L,\varphi)\in\CP}\mathsf{e}_{L,\varphi}\big\}\mpoint$$
Then clearly $\CP\subseteq \alpha\circ\beta(\CP)$. Conversely, if $\mathsf{e}_{M,\psi}\subseteq \sum_{(L,\varphi)\in\CP}\limits\mathsf{e}_{L,\varphi}$, then by Lemma~\ref{premier} there exists $(L,\varphi)\in \CP$ such that $\mathsf{e}_{M,\psi}\subseteq \mathsf{e}_{L,\varphi}$. Then $(L,\varphi)$ is a $B_K$-group, and by Corollary~\ref{inclusion BK}, this implies $(M,\psi)\twoheadrightarrow (L,\varphi)$. Hence $(M,\psi)\in \CP$, since $\CP$ is closed. Thus $\alpha\circ\beta(\CP)\subseteq \CP$, proving that $\alpha\circ\beta$ is the identity map of $\mathcal{C}l_{\BgrK{K}}$. The last assertion follows from the fact that $\mathcal{C}l_{\BgrK{K}}$ is clearly completely distributive, since its join and meet operation are union and intersection of closed subsets, respectively, and since arbitrary unions (resp. intersections) distribute over arbitrary intersections (resp. unions).\endpf
\section{Some simple $\FB_K$-modules}
\begin{mth}{Theorem} \label{SLphi}\begin{enumerate}
\item Let $(L,\varphi)$ be a $B_K$-group. Then $\mathsf{e}_{L,\varphi}$ admits a unique maximal proper subideal $\mathsf{e}_{L,\varphi}^0$, defined by
$$\mathsf{e}_{L,\varphi}^0=\sum_{\substack{(M,\psi)\in\BgrK{K}\\(M,\psi)\twoheadrightarrow (L,\varphi)\\\substack{(M,\psi)\ncong (L,\varphi)}}}\mathsf{e}_{M,\psi}\mpoint$$
\item The quotient $S_{L,\varphi}=\mathsf{e}_{L,\varphi}/\mathsf{e}_{L,\varphi}^0$ is a simple $\FB_K$-module. 
\item For any finite group $G$, let $A_G$ be a set of representatives of conjugacy classes of subgroups $X$ of $G\times K$ such that $\beta_K(X,p_2)\cong (L,\varphi)$. Then the set $\{e_X^{G\times K}\mid X\in A_G\}$ maps to an $\F$-basis of $S_{L,\varphi}(G)$ under the projection map $\mathsf{e}_{L,\varphi}(G)\to S_{L,\varphi}(G)$.
\item If $I'\subset I$ are ideals of $\FB_K$ such that $I/I'$ is a simple $\FB_K$-module, then there exists a $B_K$-group $(L,\varphi)$ such that $I/I'\cong S_{L,\varphi}$.
\end{enumerate}
\end{mth}
\pf 1. Without loss of generality, we can assume that $(L,\varphi)\in\BgrK{K}$. Using Theorem~\ref{treillis}, saying that $\mathsf{e}_{L,\varphi}$ admits a unique maximal proper subideal is equivalent to saying that the closed subset $\CP_{\mathsf{e}_{L,\varphi}}$ contains a unique maximal proper closed subset. But  
$$\CP_{\mathsf{e}_{L,\varphi}}=\{(M,\psi)\in\BgrK{K}\mid (M,\psi)\twoheadrightarrow (L,\varphi)\}\mvirg$$
so $\CP^0=\CP_{\mathsf{e}_{L,\varphi}}-\{(L,\varphi)\}$ is the unique maximal proper closed subset of $\CP_{\mathsf{e}_{L,\varphi}}$. It follows that $I_{\CP^0}= \mathsf{e}_{L,\varphi}^0$ is the unique maximal proper subideal of $\mathsf{e}_{L,\varphi}$.\mpn
2. This is clear, from 1.\mpn
3. We know from Remark~\ref{base eLphi} that $\mathsf{e}_{L,\varphi}(G)$ has a basis consisting of the idempotents $e_X^{G\times K}$, for $X$ in a set of representatives of conjugacy classes of subgroups of $G\times K$ such that $(X,p_2)\twoheadrightarrow (L,\varphi)$,
or equivalently, by Corollary~\ref{betaK}, such that $\beta_K(X,p_2)\twoheadrightarrow (L,\varphi)$. Now saying that $e_X^{G\times K}\in\mathsf{e}_{L,\varphi}^0(G)$ amounts to saying that $e_{X_{p_2}}^{X\times K}\in \mathsf{e}_{L,\varphi}^0(X)$, by Theorem~\ref{reduction1}, i.e. that $\mathsf{e}_{X,p_2}\subseteq \mathsf{e}_{M,\psi}$ for some $(M,\psi)\in\BgrK{K}$ such that $(M,\psi)\twoheadrightarrow (L,\varphi)$, but $(M,\psi)\ncong (L,\varphi)$. This in turn is equivalent to saying that $\beta_K(X,p_2)\twoheadrightarrow (L,\varphi)$, but $\beta_K(X,p_2)\ncong (L,\varphi)$. Hence $S_{L,\varphi}(G)$ has a basis consisting of the idempotents $e_X^{G\times K}$, for $X$ in a set of representatives of conjugacy classes of subgroups of $G\times K$ such that $\beta_K(X,p_2)\cong (L,\varphi)$. Assertion 2 follows.\mpn
4. Let $I'\subset I$ be ideals of $\FB_K$ such that $S=I/I'$ is a simple $\FB_K$-module, or equivalently, such that $I'$ is a maximal subideal of $I$. Then there exists $(L,\varphi)\in\BgrK{K}$ such that $\mathsf{e}_{L,\varphi}\subseteq I$ but $\mathsf{e}_{L,\varphi}\nsubseteq I'$. Hence $\mathsf{e}_{L,\varphi}+I'=I$, and $S=I/I'\cong \mathsf{e}_{L,\varphi}/(\mathsf{e}_{L,\varphi}\cap I')$. Then $\mathsf{e}_{L,\varphi}\cap I'$ is a proper subideal of $\mathsf{e}_{L,\varphi}$, so $\mathsf{e}_{L,\varphi}\cap I'\subseteq \mathsf{e}_{L,\varphi}^0$, and then $S$ maps surjectively onto $\mathsf{e}_{L,\varphi}/\mathsf{e}_{L,\varphi}^0=S_{L,\varphi}$. Since $S$ and $S_{L,\varphi}$ both are simple $\FB_K$-modules, the surjection $S\to S_{L,\varphi}$ is an isomorphism.\endpf
\begin{rem}{Remark} By Corollary~\ref{betaK mGN}, the condition $\beta_K(X,p_2)\cong (L,\varphi)$ in Assertion 3 is equivalent to the existence of a surjective morphism $s$ from $(X,p_2)$ to $(L,\varphi)$ such that $m_{X,\Ker\,s}\neq 0$. By Theorem~5.4.11 of~\cite{bisetfunctorsMSC}, or by Corollary~\ref{betaK mGN} applied to the case $K=\un$, this is equivalent to the condition $\beta(X)\cong\beta(L)$.
\end{rem} 
\begin{mth}{Corollary} Let $(L,\varphi)$ and $(M,\psi)$ be $B_K$-groups. Then the simple $\FB_K$-modules $S_{L,\varphi}$ and $S_{M,\psi}$ are isomorphic if and only if $(L,\varphi)$ and $(M,\psi)$ are isomorphic in $\grpover{K}$.
\end{mth}
\pf Clearly if $(L,\varphi)\cong(M,\psi)$ in $\grpover{K}$, then $S_{L,\varphi}\cong S_{M,\psi}$. Conversely, if $\theta:S_{L,\varphi}\to S_{M,\psi}$ is an isomorphism of $\FB_K$-modules, then for any finite group~$G$, we get an isomorphism $\theta_G:S_{L,\varphi}(G)\to S_{M,\psi}(G)$ of $\FB_K(G)$-modules. Choose $G$ such that $S_{L,\varphi}(G)\neq 0$ (e.g. $G=L$), and a subgroup $X$ of $G\times K$ such that $\beta_K(X,p_2)\cong (L,\varphi)$. Then the image $u$ of $a=e_X^{G\times K}\in \FB_K(G)$ in $S_{L,\varphi}(G)$ is non zero, and moreover $a\cdot u=u$. It follows that $\theta_G(a\cdot u)=a\cdot \theta_G(u)=\theta_G(u)$ is also non zero in $S_{M,\psi}(G)$. So there is a subgroup $Y\leq G\times K$ with $\beta_K(Y,p_2)\cong (M,\psi)$, such that the image $v$ of $e_Y^{G\times K}$ in $S_{M,\psi}(G)$ satisfies $a\cdot v\neq 0$. This forces $X$ and $Y$ to be conjugate in $G\times K$, so $(L,\varphi)\cong \beta_K(X,p_2)\cong\beta_K(Y,p_2)\cong (M,\psi)$ in $\grpover{K}$, as was to be shown.\endpf
Recall that a {\em minimal group} for a (non zero) biset functor $F$ is a finite group $G$ of minimal order such that $F(G)\neq \zero$.
\begin{mth}{Lemma} \label{SLphi non zero}Let $(L,\varphi)$ be a group over $K$.
\begin{enumerate}
\item If $N\normal L$, and $N\cap\Ker\,\varphi=\un$, then $\mathsf{e}_{L,\varphi}(L/N)\neq \zero$. 
\item If moreover $(L,\varphi)$ is a $B_K$-group, then $S_{L,\varphi}(L/N)\neq \zero$. 
\end{enumerate}
\end{mth}
\pf Indeed the map 
$$\theta:l\in L\mapsto \big(lN,\varphi(l)\big)\in (L/N)\times K$$
is injective. Let $\sur{L}\leq (L/N)\times K$ denote the image of $\theta$. Then we have a commutative diagram
$$\xymatrix{
&\sur{L}\ar[ld]_-{p_1}\ar[d]^-{p_2}\ar@{->>}[r]^-t&L\ar[d]^-\varphi\\
L/N&K\ar[r]^-i&K\makebox[0pt]{\mvirg}
}
$$
where $t:\sur{L}\to L$ is the inverse of the isomorphism $L\to \sur{L}$ induced by~$\theta$, and $i$ is the identity map of $K$. Hence $(\sur{L},p_2)\cong (L,\varphi)$ in $\grpover{K}$, and $e_{L,\varphi}=e_{\sur{L},p_2}$. In particular $e_{\sur{L}_{p_2}}^{(L/N)\times K}\in e_{L,\varphi}(L/N)$ by Theorem~\ref{reduction1}, hence $\mathsf{e}_{L,\varphi}(L/N)\neq \zero$. This proves 1. \par
If moreover $(L,\varphi)$ is a $B_K$-group, then $\beta_K(\sur{L},p_2)\cong (L,\varphi)$. It follows from Theorem~\ref{SLphi} that $e_{\sur{L}}^{(L/N)\times K}\in \mathsf{e}_{L,\varphi}(L/N)$ maps to an element of a basis of $S_{L,\varphi}(L/N)$, so $S_{L,\varphi}(L/N)\neq \zero$, proving 2.\endpf

\begin{mth}{Theorem} \label{minimal}Let $(L,\varphi)$ be a $B_K$-group, and $G$ be a finite group. The following are equivalent:
\begin{enumerate}
\item The group $G$ is a minimal group for $S_{L,\varphi}$.
\item The group $G$ is isomorphic to $L/N$, where $N$ is a normal subgroup of~$L$ of maximal order such that $N\cap\Ker\,\varphi=\un$. 
\end{enumerate}
Moreover in this case, the images in $S_{L,\varphi}(G)$ of the idempotents $e_X^{G\times K}$, where $X$ runs through a set of representatives of conjugacy classes of subgroups of $G\times K$ such that $(X,p_2)\cong (L,\varphi)$, form an $\F$-basis of $S_{L,\varphi}(G)$.
\end{mth}
\pf By Theorem~\ref{SLphi}, saying that $S_{L,\varphi}(G)\neq\zero$ for a finite group~$G$ amounts to saying that there exists a subgroup $X$ of $G\times K$ such that $\beta_K(X,p_2)\cong (L,\varphi)$ in $\grpover{K}$. Equivalently, there is a commutative diagram
\begin{equation}\label{diagram}
\vcenter{
\xymatrix{
&X\ar[ld]_-{p_1}\ar[d]^-{p_2}\ar@{->>}[r]^-s&L\ar[d]^-\varphi\\
G&K\ar[r]^-i&K\makebox[0pt]{\mvirg}
}
}
\end{equation}
where \begin{itemize}
\item $s$ is surjective and $m_{X,\Ker\,s\neq 0}$,
\item $i$ is an inner automorphism of $K$,
\item the map $(p_1,p_2):X\to (G\times K)$ is injective.
\end{itemize}
Now we proceed with the proof of Theorem~\ref{minimal}.\par\noindent
\fbox {$1\Rightarrow 2$} If $G$ is minimal for $S_{L,\varphi}$, then $S_{L,\varphi}(G)\neq \zero$, so we have a diagram~(\ref{diagram}). Let $H=p_1(G)$. Replacing $G$ by $H$ in this diagram gives a diagram for the group $H$ with the same properties, so $S_{L,\varphi}(H)\neq 0$. Hence $H=G$ by minimality of~$G$. In other words $p_1$ is surjective, so $G\cong X/\Ker\,p_1$. Let $N=s(\Ker\,p_1)$. If $u\in N\cap\Ker\,\varphi$, then $u=s(x)$ for some $x\in X$, and then $\varphi\circ s(x)=i\circ p_2(x)=1$, so $p_2(x)=1$. Thus $x=1$ since $\Ker\,p_1\cap\Ker\,p_2=\un$. Moreover $N$ in normal in $L$, since $s$ is surjective. Lemma~\ref{SLphi non zero} shows that $S_{L,\varphi}(L/N)\neq \zero$, and by minimality of $G$, the surjection $\sur{s}:G\cong X/\Ker\,p_1\twoheadrightarrow L/N$ induced by $s$ must be an isomorphism. Lemma~\ref{SLphi non zero} also implies that $N$ is a normal subgroup of maximal order of $L$ such that $N\cap\Ker\,\varphi$. Hence 2 holds.\par 
Observe that it also follows that $\Ker\,s\leq\Ker\,p_1$, so $\Ker\,s=\un$ since $\Ker\,s\leq\Ker\,p_2$ as $\varphi\circ s=i\circ p_2$, and $\Ker\,p_1\cap\Ker\,p_2=\un$. So $s$ is an isomorphism $X\to L$. This proves the last assertion of the theorem.\mpn
\fbox{$2\Rightarrow 1$} Suppose that 2 holds. Then $S_{L,\varphi}(G)\neq 0$, by the above claim. By the first part of the proof, if $H$ is a minimal group for $S_{L,\varphi}$, then $H\cong L/M$, where $M$ is a normal subgroup of maximal order such that $M\cap \Ker\,\varphi=\un$. Then $|M|=|N|$, so $|G|=|H|$, and $S_{L,\varphi}(G')=\zero$ for any group $G'$ of order smaller than $|G|=|H|$. Hence $G$ is minimal for $S_{L,\varphi}$, and 1 holds.\endpf 
\pagebreak[3]
\begin{mth}{Corollary} Let $(L,\varphi)$ be a group over $K$. The following conditions are equivalent:
\begin{enumerate}
\item $\varphi:L\to K$ is injective.
\item $(L,\varphi)$ is a $B_K$-group and $S_{L,\varphi}(\un)\neq \zero$.
\end{enumerate}
\end{mth}
\pf \fbox{$1\Rightarrow 2$} If $\varphi$ is injective, then $(L,\varphi)$ is a $B_K$-group (cf. Example~\ref{exBK}). Moreover $L\cap\Ker\,\varphi=\un$, so $S_{L,\varphi}(L/L)= S_{L,\varphi}(\un)\neq \zero$.\mpn
\fbox{$2\Rightarrow 1$} If $(L,\varphi)$ is a $B_K$-group and $S_{L,\varphi}(\un)\neq \zero$, then $\un$ is a minimal group for $S_{L,\varphi}$. So there is a normal subgroup $N$ of $L$ of maximal order such that $N\cap\Ker\,\varphi=\un$, such that moreover $L/N\cong \un$. Hence $N=L$, and $\Ker\,\varphi=N\cap\Ker\,\varphi=\un$.\endpf
\begin{rem}{Example} Let $L=C_2\times (C_3\rtimes C_4)$ be a direct product of a group of order 2, generated by the element $a$, and a semidirect product of a group of order 3, generated by $b$, and a cyclic group of order 4, generated by $c$ (so $cbc^{-1}=b^{-1}$). Let $P$ be the subgroup of $L$ generated by $a$ and $b$. Then $P$ is cyclic of order 6, and the factor group $K=L/P$ is cyclic of order 4, generated by the class $cP$. Let $\varphi:L\to K$ be the projection map. One can check that $(L,\varphi)$ is a $B_K$-group, i.e. that $m_{L,Q}=0$ when $Q$ is any of the non trivial subgroups of $P$ (these subgroups are all normal in $L$, as $P$ is cyclic).\par
Then the subgroups $M=\langle ac^2\rangle$ and $N=\langle c^2\rangle$ both are normal (central, in fact) subgroups of $L$ of maximal order (equal to 2) intersecting trivially $P=\Ker\,\varphi$. So the groups $G=L/M$ and $H=L/N$ are both minimal groups\footnote{One can show moreover that $S_{L,\varphi}(G)$ and $S_{L,\varphi}(H)$ are both one dimensional.} for the simple $\FB_K$-module $S_{L,\varphi}$, but they are not isomorphic, as $G\cong C_3\rtimes C_4$ but $H\cong C_2\times S_3$, where $S_3$ is the symmetric group of degree~3. This gives yet another counterexample to a conjecture I made in 2010, saying that the minimal groups for a Green biset functor should form a single isomorphism class of groups. The first counterexample to this conjecture was found by Nadia Romero in 2013 (cf.~\cite{romero-fibred}). Another counterexample was found recently by Ibrahima Tounkara (cf. \cite{tounkara-jofa}).
\end{rem}
\section{Restriction to $p$-groups}
In this section, we fix a prime number $p$, and restrict the functor $\FB_K$ to finite $p$-groups. We obtain a Green $p$-biset functor $\FB_{K}^{(p)}$. We do not assume that $K$ is itself a $p$-group. \par
In order to study the ideals of $\FB_K^{(p)}$, it is natural to try to determine those groups $(L,\varphi)$ over $K$ for which the restriction of $\mathsf{e}_{L,\varphi}$ to $p$-groups does not vanish. This motivates the following definition:
\begin{mth}{Definition} Let $K$ be a finite group. Then a group $(L,\varphi)$ over $K$ is called {\em $p$-persistent} if there is a finite $p$-group $P$ such that $\mathsf{e}_{L,\varphi}(P)\neq\zero$.\par
We denote by $\grpover{K}^{(p)}$ the full subcategory of $\grpover{K}$ consisting of $p$-persistent groups over $K$.

\end{mth}
\begin{rem}{Remarks} \label{quotient pepere}\begin{enumerate}
\item If $X$ is a subgroup of $P\times K$, where $P$ is a $p$-group, then $(X,p_2)$ is $p$-persistent: indeed $e_X^{P\times K}\in\mathsf{e}_{X,p_2}(P)$ by Corollary~\ref{reduction}.
\item Any quotient of a $p$-persistent group over $K$ is $p$-persistent: indeed is $s:(M,\psi)\twoheadrightarrow (L,\varphi)$ is a surjective morphism in $\grpover{K}$, then $\mathsf{e}_{M,\psi}\subseteq \mathsf{e}_{L,\varphi}$ by Lemma~\ref{quotient}. It follows that $\mathsf{e}_{L,\varphi}(P)\neq \zero$ if $P$ is a $p$-group such that $\mathsf{e}_{M,\psi}(P)\neq \zero$. In particular, if $(L,\varphi)$ is $p$-persistent, then $\beta_K(L,\varphi)$ is a $p$-persistent $B_K$-group.
\end{enumerate}
\end{rem}
\begin{mth}{Notation} When $L$ is a finite group, we denote by $O^p(L)$ the subgroup of $L$ generated by $p'$-elements, and by $L^{[p]}$ the quotient $L/O^p(L)$.
\end{mth}
Recall that $O^p(L)$ is the smallest normal subgroup $N$ of $L$ such that $L/N$ is a $p$-group. Also recall that if $s:M\to L$ is a surjective group homomorphism, then $s\big(O^p(M)\big)=O^p(L)$. Indeed $N=s\big(O^p(M)\big)\normal L$, and $s$ induces a surjection $M^{[p]}\to L/N$. So $L/N$ is a $p$-group, thus $N\geq O^p(L)$. But $N$ is generated by $p'$-elements, as $O^p(M)$ is, so $N\leq O^p(L)$.
\begin{mth}{Proposition} \label{persistent}Let $(L,\varphi)$ be a group over $K$. The following are equivalent:
\begin{enumerate}
\item $(L,\varphi)$ is $p$-persistent.
\item $\mathsf{e}_{L,\varphi}\big(L^{[p]}\big)\neq\zero$.
\item $m_{L,O^p(L)\cap\Ker\,\varphi}\neq 0$.
\end{enumerate}
\end{mth}
\pf Indeed if 3 holds, then setting $N=O^p(L)\cap \Ker\,\varphi$, we have $\mathsf{e}_{L,\varphi}=\mathsf{e}_{L/N,\varphi/N}$ by Theorem~\ref{surjection forte}. Moreover $O^p(L/N)=O^p(L)/N$, and $\Ker\,(\varphi/N)=\Ker\,\varphi/N$. Thus $O^p(L/N)\cap \Ker(\varphi/N)=\un$, so $\mathsf{e}_{L/N,\varphi/N}\big((L/N)/O^p(L/N)\big)$ is non zero by Lemma~\ref{SLphi non zero}. But 
$$\mathsf{e}_{L/N,\varphi/N}\big((L/N)/O^p(L/N)\big)\cong \mathsf{e}_{L/N,\varphi/N}\big(L/O^p(L)\big)=\mathsf{e}_{L,\varphi}(L^{[p]})\mvirg$$
so 2 holds. Clearly 2 implies 1, as $L^{[p]}$ is a $p$-group. Now if 1 holds, let $P$ be a $p$-group such that $\mathsf{e}_{L,\varphi}(P)\neq\zero$. Let $N$ be a normal subgroup of $L$ contained in $\Ker\,\varphi$, and maximal such that $m_{L,N}\neq 0$. Then setting $\sur{L}=L/N$ and $\sur{\varphi}=\varphi/N$, we have $\beta_K(L,\varphi)\cong (\sur{L},\sur{\varphi})$, and $\mathsf{e}_{L,\varphi}=\mathsf{e}_{\sur{L},\sur{\varphi}}$ by Theorem~\ref{surjection forte}. Moreover as $(\sur{L},\sur{\varphi})$ is a $B_K$-group, by Theorem~\ref{eLphi}, there exists a subgroup $X$ of $P\times K$, and a commutative diagram
\begin{equation*}\vcenter{
\xymatrix{
&X\ar[ld]_-{p_1}\ar[d]^-{p_2}\ar@{->>}[r]^-s&\sur{L}\ar[d]^-{\sur{\varphi}}\\
P&K\ar[r]^-i&K\makebox[0pt]{\mvirg}
}
}
\end{equation*}
where $s$ is surjective and $i$ is an inner automorphism of $K$. Then $N=s(\Ker\,p_1)$ is a normal subgroup of $\sur{L}$, as $s$ is surjective. Moreover if $l\in N\cap\Ker\,\sur{\varphi}$, then $l=s(x)$ for some $x\in\Ker\,p_1$, so $p_1(x)=1$ and $i\circ p_2(x)=\varphi\circ s(x)=1$, so $p_2(x)=1$. Hence $x=1$, and $l=1$, so $N\cap \Ker\,\sur{\varphi}=1$. Now $s$ induces a surjection $X/\Ker\,p_1\cong p_1(X)\twoheadrightarrow L/N$, so $L/N$ is a $p$-group, thus $N\geq O^p(L)$. It follows that $O^p(\sur{L})\cap\Ker\,\sur{\varphi}=\un$. Now if $\pi:L\to \sur{L}=L/N$ is the projection map, we have $\sur{\varphi}\circ\pi=\varphi$, so
$$\pi\big(O^p(L)\cap\Ker\,\varphi\big)\leq O^p(\sur{L})\cap\Ker\,\sur{\varphi}=\un\mvirg$$
that is $O^p(L)\cap\Ker\,\varphi\leq N=\Ker\,\pi$. Then if $M=O^p(L)\cap\Ker\,\varphi$, we have $m_{L,N}=m_{L,M}m_{L/M,N/M}\neq 0$, hence $m_{L,M}\neq 0$, so 3 holds.
\endpf
\begin{mth}{Corollary} \label{condition BK}Let $(L,\varphi)$ be a $p$-persistent $B_K$-group. Then 
$$O^p(L)\cap\Ker\,\varphi=\un\mpoint$$
\end{mth}
\pf Indeed $m_{L,O^p(L)\cap\Ker\,\varphi}\neq 0$, and $(L,\varphi)$ is a $B_K$-group.\endpf
\begin{mth}{Notation} When $(L,\varphi)$ is a $p$-persistent group over $K$, we denote by $L^{(p)}_\varphi$ the subgroup of $L^{[p]}\times K$ defined by
$$ L^{(p)}_\varphi=\big\{\big(lO^p(L),\varphi(l)\big)\mid l\in L\big\}\mpoint$$
\end{mth}
The following theorem is analogous to Theorem~\ref{reduction1}:
\pagebreak[3]
\begin{mth}{Theorem} \label{reduction1p}Let $I$ be an ideal of the Green biset functor $\FB_K^{(p)}$.
 If $G$ is a finite $p$-group and $L$ is a subgroup of $G\times K$, the following conditions are equivalent:
\begin{enumerate}
\item The idempotent $e_{L}^{G\times K}$ belongs to $I(G)$.
\item The idempotent $e_{L_{p_2}^{(p)}}^{L^{[p]}\times K}$ belongs to $I(L^{[p]})$.
\end{enumerate}
\end{mth}
\pf The proof is similar to the proof of Theorem~\ref{reduction1}, so we only sketch~it. If $L\leq G\times K$, denote by $\widehat{L}$ the image of $L$ in the group $L^{[p]}\times G$ by the map $l\mapsto \big(lO^p(L),p_1(l)\big)$. Recall that $\Ker\,p_1\geq O^p(L)$, since $G$ is a $p$-group. Furthermore $p_1(\widehat{L})=L^{[p]}$, $k_1(\widehat{L})=\Ker\,p_1/O^p(L)$, $p_2(\widehat{L})=p_1(L)$, and $k_2(\widehat{L})=p_1\big(O^p(L)\big)=\un$. The $(L^{[p]},G)$-biset $U=(L^{[p]}\times G)/\widehat{L}$ factors as
$$U\cong \Inf_{L^{[p]}/k_1(\widehat{L})}^{L^{[p]}}\circ\Iso(\theta^{-1})\circ\Res_{p_1(L)}^G\mvirg$$
where $\theta:L^{[p]}/k_1(\widehat{L})\to p_1(G)$ is the isomorphism induced by the map $l O^p(L)\mapsto p_1(l)$ from $L^{[p]}$ to $G$.\par
If $e_{L}^{G\times K}$ belongs to $I(G)$, then $\FB_K^{(p)}(U)(e_{L}^{G\times K})$ belongs to $I(L^{[p]})$. As in the proof of Theorem~\ref{reduction1}, one can check that the product $e_{L_{p_2}^{(p)}}^{L^{[p]}\times K}\cdot \FB_K^{(p)}(e_{L}^{G\times K})$ is non zero. As it is a scalar multiple of $e_{L_{p_2}^{(p)}}^{L^{[p]}\times K}$, we get that $e_{L_{p_2}^{(p)}}^{L^{[p]}\times K}\in I(L^{[p]})$, thus 1 implies 2.\par
Conversely, assume that $e_{L_{p_2}^{(p)}}^{L^{[p]}\times K}\in I(L^{[p]})$. Then, as in the proof of Theorem~\ref{reduction1} again, the opposite biset $U\op$ factors as
$$U\op\cong \Ind_{p_1(L)}^G\circ \Iso(\theta)\circ \Def_{L^{[p]}/k_1(\widehat{L})}^{L^{[p]}}\mvirg$$
and the element $\FB_K^{(p)}(U\op)\big(e_{L_{p_2}^{(p)}}^{L^{[p]}\times K}\big)$ belongs to $I(G)$. One can can check moreover that there is a non zero scalar $\lambda$ such that
$$\FB_K^{(p)}(U\op)\big(e_{L_{p_2}^{(p)}}^{L^{[p]}\times K}\big)=\lambda\,m_{L_{p_2}^{(p)},L_{p_2}^{(p)}\cap(N\times \un)}\,e_L^{G\times K}\mvirg$$
where $N=k_1(\widehat{L})=\Ker\,p_1/O^p(L)\leq L^{[p]}$. \par
But if $\big(lO^p(L),p_2(l)\big)\in L_{p_2}^{(p)}\cap(N\times \un)$, then $l\in\Ker\,p_2\cap\Ker\,p_1=\un$. It follows that $m_{L_{p_2}^{(p)},L_{p_2}^{(p)}\cap(N\times \un)}=m_{L_{p_2}^{(p)},\un}=1$, and $e_L^{G\times K}\in I(G)$, as $\lambda\neq 0$. Hence 2 implies 1.\endpf
\pagebreak[3]
\begin{mth}{Corollary} \label{reductionp} Let $G$ be a finite $p$-group, and $L$ be a subgroup of $G\times\nolinebreak K$. Then the ideal of $\FB_K^{(p)}$ generated by $e_L^{G\times K}$ is equal to the ideal of $\FB_K^{(p)}$ generated by $e_{L_{p_2}^{(p)}}^{L^{[p]}\times K}$
\end{mth}
\pf The proof is the same as the proof of Corollary~\ref{reduction}.\endpf
\begin{mth}{Notation} Let $(L,\varphi)$ be a $p$-persistent group over $K$. We denote by $\mathsf{e}_{L,\varphi}^{(p)}$ the ideal of $\FB_K^{(p)}$ generated by $e_{L_\varphi^{(p)}}^{L^{[p]}\times K}\in\FB_K^{(p)}(L^{[p]})$.
\end{mth}
\begin{mth}{Theorem} \label{betaKp}Let $s:(M,\psi)\twoheadrightarrow (L,\varphi)$ be a surjective morphism in $\grpover{K}$, and assume that $(M,\psi)$ is $p$-persistent. Then:
\begin{enumerate}
\item $(L,\varphi)$ is $p$-persistent, and $\mathsf{e}_{M,\psi}^{(p)}\subseteq \mathsf{e}_{L,\varphi}^{(p)}$.
\item If $m_{M,\Ker\,s}\neq 0$, then $\mathsf{e}_{M,\psi}^{(p)}= \mathsf{e}_{L,\varphi}^{(p)}$.
\end{enumerate}
\end{mth}
\pf 1. We already observed in Remarks~\ref{quotient pepere} that any quotient of a $p$-persistent group over $K$ is itself $p$-persistent, hence $(L,\varphi)$ is $p$-persistent. Let $i$ be an inner automorphism of $K$ such that $i\circ\psi=\varphi\circ s$. The surjection $s:M\to L$ induces a surjection $s^{[p]}:M^{[p]}\to L^{[p]}$, hence a surjection 
$$s^{[p]}\times\Id_K:M^{[p]}\times K\to L^{[p]}\times K\mpoint$$
Let $u=\big(mO^p(M),\psi(m)\big)$ be the image of $m\in M$ in $M_\psi^{(p)}$. Then
$$(s^{[p]}\times\Id_K)(u)=\big(s(m)O^p(L),\psi(m)\big)=\Big(s(m)O^p(L),i^{-1}\circ\varphi\big(s(m)\big)\Big)\mvirg$$
which shows that $s^{[p]}\times\Id_K$ maps $M_\psi^{(p)}$ to a conjugate of $L_\varphi^{(p)}$ %
 in $L^{[p]}\times K$. Then the idempotent $e_{M_\psi^{(p)}}^{M^{[p]}\times K}$ appears in the decomposition of 
$$\FB_K^{(p)}\big(\Inf_{M^{[p]}/\Ker\,s^{[p]}}^{M^{[p]}}\circ\Iso(\alpha^{-1})\big)(e_{L_\varphi^{(p)}}^{L^{[p]}\times K})\mvirg$$
where $\alpha:M^{[p]}/\Ker\,s^{[p]}\to L^{[p]}$ is the canonical isomorphism. It follows that $e_{M_\psi^{(p)}}^{M^{[p]}\times K}\in\mathsf{e}_{L,\varphi}^{(p)}(M^{[p]})$, hence $\mathsf{e}_{M,\psi}^{(p)}\subseteq\mathsf{e}_{L,\varphi}^{(p)}$.\mpn
2. Consider now $v=\FB_K^{(p)}\big(\Iso(\alpha)\circ\Def_{M^{[p]}/\Ker\,s^{[p]}}^{M^{[p]}}\big)(e_{M_\psi^{(p)}}^{M^{[p]}\times K})\in \mathsf{e}_{M,\psi}^{(p)}(L^{[p]})$. By Lemma~\ref{deflation}, there is a non zero scalar $\lambda$ such that
\begin{equation}\label{v}
v=\lambda\,m_{M_\psi^{(p)},M_\psi^{(p)}\cap(\Ker\,s^{[p]}\times\un)}e_{L_\varphi^{(p)}}^{L^{[p]}\times K}\mpoint
\end{equation}
Now the projection $m\in M\mapsto \big(mO^p(M),\psi(m)\big)\in M_\psi^{(p)}$ induces an isomorphism $M_\psi^{(p)}\cong M/\big(O^p(M)\cap\Ker\,\psi\big)$. As $\Ker\,s^{[p]}=\Ker\,s\,O^p(L)/O^p(L)$, the subgroup $M_\psi^{(p)}\cap(\Ker\,s^{[p]}\times\un)$ maps to $(\Ker\,s\,O^p(M)\cap\Ker\,\psi)/\big(O^p(M)\cap\Ker\,\psi\big)$ under this isomorphism. \par
Moreover $\Ker\,s\,O^p(M)\cap\Ker\,\psi=\Ker\,s\big(O^p(M)\cap\Ker\,\psi\big)$ as $\Ker\,s\leq\Ker\,\psi$. It follows that
$$m_{M_\psi^{(p)},M_\psi^{(p)}\cap(\Ker\,s^{[p]}\times\un)}=m_{M/(O^p(M)\cap\Ker\,\psi),\Ker\,s(O^p(M)\cap\Ker\,\psi)/(O^p(M)\cap\Ker\,\psi)}\mpoint$$
Multiplying by $m_{M,O^p(M)\cap\Ker\,\psi}$, which is non zero by Proposition~\ref{persistent}, since $(M,\psi)$ is $p$-persistent, this gives
\begin{align*}
m_{M,O^p(M)\cap\Ker\,\psi}\,m_{M_\psi^{(p)},M_\psi^{(p)}\cap(\Ker\,s^{[p]}\times\un)}&=m_{M,\Ker\,s(O^p(M)\cap\Ker\,\psi)}\\
&=m_{M,\Ker\,s}\,m_{M/\Ker\,s,\Ker\,s(O^p(M)\cap\Ker\,\psi)/\Ker\,s}\\
&=m_{M,\Ker\,s}\,m_{L,O^p(L)\cap\Ker\varphi}\mvirg
\end{align*}
as the canonical isomorphism $M/\Ker\,s\to L$ maps $\Ker\,s(O^p(M)\cap\Ker\,\psi)/\Ker\,s$ to $O^p(L)\cap\Ker\varphi$. Since $m_{L,O^p(L)\cap\Ker\varphi}\neq 0$ as $(L,\varphi)$ is $p$-persistent, and since $m_{M,\Ker\,s}\neq 0$ by assumption, it follows that $m_{M_\psi^{(p)},M_\psi^{(p)}\cap(\Ker\,s^{[p]}\times\un)}\neq 0$, hence $e_{L_\varphi^{(p)}}^{L^{[p]}\times K}$ is a non zero scalar multiple of $v$, by~\ref{v}. It follows that $e_{L_\varphi^{(p)}}^{L^{[p]}\times K}$ belongs to $\mathsf{e}_{M,\psi}^{(p)}(L^{[p]})$, so $\mathsf{e}_{L,\varphi}^{(p)}\subseteq \mathsf{e}_{M,\psi}^{(p)}$, and $\mathsf{e}_{L,\varphi}^{(p)}= \mathsf{e}_{M,\psi}^{(p)}$, as was to be shown.\endpf
\begin{mth}{Corollary} \label{restriction eLphi}Let $(L,\varphi)$ be a $p$-persistent group over $K$. Then the restriction of $\mathsf{e}_{L,\varphi}$ to finite $p$-groups is equal to $\mathsf{e}_{L,\varphi}^{(p)}$.
\end{mth}
\pf Since $\mathsf{e}_{L,\varphi}=\mathsf{e}_{\beta_K(L,\varphi)}$  by Corollary~\ref{betaK}, and since $\mathsf{e}_{L,\varphi}^{(p)}=\mathsf{e}_{\beta_K(L,\varphi)}^{(p)}$ by Theorem~\ref{betaKp}, we may assume that $(L,\varphi)$ is a $B_K$-group. By Corollary~\ref{condition BK}, we have $O^p(L)\cap\Ker\,\varphi=\un$. Thus the projection $L\to L_\varphi^{(p)}$ is an isomorphism, and it induces an isomorphism $(L,\varphi)\cong (L_\varphi^{(p)},p_2)$. Hence $e_{L_\varphi^{(p)}}^{L^{(p]}\times K}\in \mathsf{e}_{L,\varphi}(L^{[p]})$, and $\mathsf{e}_{L,\varphi}^{(p)}$ is contained in the restriction of $\mathsf{e}_{L,\varphi}$ to $p$-groups.\par
Conversely, if $G$ is a $p$-group and $e_X^{G\times K}\in \mathsf{e}_{L,\varphi}(G)$, then $(X,p_2)\twoheadrightarrow (L,\varphi)$ by Theorem~\ref{eLphi}. Then $\mathsf{e}_{X,p_2}^{(p)}\subseteq \mathsf{e}_{L,\varphi}^{(p)}$, hence $e_{X}^{G\times K}\in \mathsf{e}_{L,\varphi}^{(p)}$ by Corollary~\ref{reductionp}. Hence the restriction of $\mathsf{e}_{L,\varphi}$ is contained in $\mathsf{e}_{L,\varphi}^{(p)}$, which completes the proof.\endpf
\begin{mth}{Corollary} \label{inclusion BKp}Let $(L,\varphi)$ be a $p$-persistent $B_K$-group, and $(M,\psi)$ be a $p$-persistent group over~$K$. Then $\mathsf{e}_{M,\psi}^{(p)}\subseteq \mathsf{e}_{L,\varphi}^{(p)}$ if and only if $(M,\psi)\twoheadrightarrow (L,\varphi)$.\vspace{-2ex}
\end{mth}
\pf Indeed if $(M,\psi)\twoheadrightarrow (L,\varphi)$, then $\mathsf{e}_{M,\psi}^{(p)}\subseteq \mathsf{e}_{L,\varphi}^{(p)}$ by Theorem~\ref{betaKp}. Conversely, if $\mathsf{e}_{M,\psi}^{(p)}\subseteq \mathsf{e}_{L,\varphi}^{(p)}$, showing that $(M,\psi)\twoheadrightarrow (L,\varphi)$ amounts to showing that $\beta_K(M,\psi)\twoheadrightarrow (L,\varphi)$, because $(L,\varphi)$ is a $B_K$-group. Now $\mathsf{e}_{M,\psi}=\mathsf{e}_{\beta_K(M,\psi)}$, hence $\mathsf{e}_{M,\psi}^{(p)}=\mathsf{e}_{\beta_K(M,\psi)}^{(p)}$ by Corollary~\ref{restriction eLphi}, and we can assume that $(M,\psi)$ is also a $B_K$-group. \par
If $\mathsf{e}_{M,\psi}^{(p)}\subseteq \mathsf{e}_{L,\varphi}^{(p)}$, then $e_{M_\psi^{(p)}}^{M^{[p]}\times K}\in \mathsf{e}_{L,\varphi}^{(p)}(M^{[p]})$, and $\mathsf{e}_{L,\varphi}^{(p)}(M^{[p]})=\mathsf{e}_{L,\varphi}(M^{[p]})$ by Corollary~\ref{restriction eLphi}. Hence $(M_\psi^{(p)},p_2)\twoheadrightarrow (L,\varphi)$ by Theorem~\ref{eLphi}. But the projection $M\to M_\psi^{(p)}$ is a group isomorphism, since $(M,\psi)$ is a $B_K$-group. It is in fact an isomorphism  from $(M,\psi)$ to $(M_\psi^{(p)},p_2)$ in $\grpover{K}$. It follows that $(M,\psi)\twoheadrightarrow (L,\varphi)$. \endpf
The following is analogous to Lemma~\ref{premier}, and the proof is the same:
\begin{mth}{Lemma} \label{premierp}Let $\mathcal{A}$ be a set of ideals of $\FB_K^{(p)}$, and $(M,\psi)$ be a $p$-persistent group over~$K$. The following are equivalent:
\begin{enumerate}
\item $\mathsf{e}_{M,\psi}^{(p)}\subseteq \sum_{I\in\mathcal{A}}\limits I$.
\item There exists $I\in\mathcal{A}$ such that $\mathsf{e}_{M,\psi}^{(p)}\subseteq I$.
\end{enumerate}
\end{mth}
\pf Clearly 2 implies 1. Now 1 is equivalent to saying that 
$$e_{M_\psi^{(p)}}^{M^{[p]}\times K}\in \sum_{I\in\mathcal{A}}\limits I(M^{[p]})\mpoint$$
If this holds, there exists $I\in\mathcal{A}$ and $u\in I(M^{[p]})$ such that $e_{M_\psi^{(p)}}^{M^{[p]}\times K}\cdot u\neq 0$. Now $e_{M_\psi^{(p)}}^{M^{[p]}\times K}\cdot u\in I(M^{[p]})$, and moreover there is a scalar $\lambda\in\F$ such that $e_{M_\psi^{(p)}}^{M^{[p]}\times K}\cdot u=\lambda e_{M_\psi^{(p)}}^{M^{[p]}\times K}\neq 0$. Hence $\lambda\neq 0$, and $e_{M_\psi^{(p)}}^{M^{[p]}\times K}\in I(M^{[p]})$. In other words $\mathsf{e}_{M,\psi}^{(p)}\subseteq I$, so 1 implies 2.\endpf 
\begin{mth}{Notation} Let $\BgrKp{K}{p}$ denote the subset of $\BgrK{K}$ consisting of $p$-persistent $B_K$-groups. 
\end{mth}
As before, a subset $\mathcal{P}$ of $\BgrKp{K}{p}$ is called {\em closed} if $$\forall (L,\varphi)\in\CP,\;\forall (M,\psi)\in \BgrKp{K}{p},\; (M,\psi)\twoheadrightarrow (L,\varphi)\;\implies \;(M,\psi)\in\CP\mpoint$$
\pagebreak[3]
\begin{mth}{Theorem} \label{treillisp}Let $\mathcal{I}_{\FB_K^{(p)}}$ be the lattice  of ideals of $\FB_K^{(p)}$, ordered by inclusion of ideals, and $\mathcal{C}l_{\BgrKp{K}{p}}$ be the lattice of closed subsets of $\BgrKp{K}{p}$, ordered by inclusion of subsets. Then the map
$$I\in\mathcal{I}_{\FB_K^{(p)}}\mapsto\CP_I=\{(L,\varphi)\in\BgrKp{K}{p}\mid \mathsf{e}_{L,\varphi}^{(p)}\subseteq I\}$$
is an isomorphism of lattices from $\mathcal{I}_{\FB_K^{(p)}}$ to $\mathcal{C}l_{\BgrKp{K}{p}}$. The inverse isomorphism is the map
$$\CP\in \mathcal{C}l_{\BgrKp{K}{p}}\mapsto I_\CP=\sum_{(L,\varphi)\in \CP}\mathsf{e}_{L,\varphi}^{(p)}\mpoint$$
In particular $\mathcal{I}_{\FB_K^{(p)}}$ is completely distributive.
\end{mth}
\pf First the map $I\in\mathcal{I}_{\FB_K^{(p)}}\mapsto \CP_I\in \mathcal{C}l_{\BgrKp{K}{p}}$ is well defined: indeed $\CP_I\in \mathcal{C}l_{\BgrKp{K}{p}}$ by Theorem~\ref{betaKp}. This map is obviously order preserving. Similarly, the map $\CP\in \mathcal{C}l_{\BgrKp{K}{p}}\mapsto I_\CP=\sum_{(L,\varphi)\in \CP}\limits\mathsf{e}_{L,\varphi}^{(p)}$ is also well defined and order preserving.\par
Hence all we need to show is that if $I$ is an ideal of $\FB_K^{(p)}$, then
\begin{equation}\label{a voir 1}I=\sum_{(L,\varphi)\in \CP_I}\mathsf{e}_{L,\varphi}^{(p)}\mvirg
\end{equation}
and that if $\CP$ is a closed subset of $\BgrKp{K}{p}$, and $(M,\psi)\in \BgrKp{K}{p}$, then 
\begin{equation}\label{a voir 2}\mathsf{e}_{M,\psi}^{(p)}\subseteq \sum_{(L,\varphi)\in \CP}\mathsf{e}_{L,\varphi}^{(p)}\;\;\Leftrightarrow\;\; (M,\psi)\in \CP\mpoint
\end{equation}
\hspace{-1ex}For~\ref{a voir 1}, let $J=\sum_{(L,\varphi)\in \CP_I}\limits\mathsf{e}_{L,\varphi}^{(p)}$. Then $J\subseteq I$ by definition of $\CP_I$. Conversely, let $G$ be a finite $p$-group, and $u=\sum_{X\in E}\limits\lambda_Xe_X^{G\times K}$ be an element of $I(G)$, where $\lambda_X\in\F$, and $E$ is a set of representatives of conjugacy classes of subgroups of $G\times K$. Then $e_X^{G\times K}\cdot u=\lambda_Xe_X^{G\times K}\in I(G)$, for any $X\in E$. So if $\lambda_X\neq 0$, then $e_X^{G\times K}\in I(G)$. Equivalently, by Theorem~\ref{reduction1p}, $e_{X_{p_2}^{(p)}}^{X^{[p]}\times K}\in I(X^{[p]})$, that is $\mathsf{e}_{X,p_2}^{(p)}\subseteq I$. Let $(L,\varphi)$ be the element of $\BgrKp{K}{p}$ isomorphic to $\beta_K(X,{p_2})$. Then $\mathsf{e}_{X,p_2}^{(p)}=\mathsf{e}_{L,\varphi}^{(p)}$ by Theorem~\ref{betaKp}, and $(L,\varphi)\in \CP_I$. \par
Moreover $e_{X_{p_2}^{(p)}}^{X^{[p]}\times K}\in \mathsf{e}_{L,\varphi}^{(p)}(X^{[p]})$, or equivalently $e_X^{G\times K}\in \mathsf{e}_{L,\varphi}^{(p)}(G)\subseteq J(G)$. As this holds for any $X\in E$ such that $\lambda_X\neq 0$, we have also $u\in J(G)$, so $J(G)=I(G)$, as $u$ was arbitrary in $I(G)$, and $J=I$, as $G$ was an arbitrary finite $p$-group. This completes the proof of~\ref{a voir 1}.\mpn
As for~\ref{a voir 2}, clearly if $(M,\psi)\in \CP$, then  $\mathsf{e}_{M,\psi}^{(p)}\subseteq \sum_{(L,\varphi)\in \CP}\limits\mathsf{e}_{L,\varphi}^{(p)}$. Conversely if $\mathsf{e}_{M,\psi}^{(p)}\subseteq \sum_{(L,\varphi)\in \CP}\limits\mathsf{e}_{L,\varphi}^{(p)}$, then by Lemma~\ref{premierp}, there exists $(L,\varphi)\in \CP$ such that $\mathsf{e}_{M,\psi}^{(p)}\subseteq \mathsf{e}_{L,\varphi}^{(p)}$. Hence $(M,\psi)\twoheadrightarrow (L,\varphi)$, by Corollary~\ref{inclusion BKp}. Since $(L,\varphi)\in \CP$ and $\CP$ is closed, we get that $(M,\psi)\in\CP$, as was to be shown.\endpf
\begin{mth}{Theorem} \label{structure BKp groups}Let $(L,\varphi)$ be a $p$-persistent $B_K$-group. Let $[s_K]$ be a set of representatives of conjugacy classes of subgroups of $K$. Let $H$ be the unique element of $[s_K]$ conjugate to $\varphi(L)$, and $j_H:H\hookrightarrow K$ be the inclusion map. Then one and one only of the following holds:
\begin{enumerate}
\item $\Ker\,\varphi=\un$, and $(L,\varphi)\cong (H,j_H)$ in $\grpover{K}$.
\item $\Ker\varphi\cong C_p$, the group $H^{[p]}$ is cyclic and non trivial, and $(L,\varphi)\cong (C_p\times H,j_H\circ\pi_H)$ in $\grpover{K}$, where $\pi_H:C_p\times H\to K$ is the projection onto $H$.
\item $\Ker\varphi\cong C_p\times C_p$, the group $H^{[p]}$ is trivial - in other words $H$ is a $p$-perfect subgroup of $K$~- and $(L,\varphi)\cong (C_p\times C_p\times H,j_H\circ\pi_H)$ in $\grpover{K}$, where $\pi_H:C_p\times C_p\times H\to K$ is the projection onto $H$.
\end{enumerate}
\end{mth}
\pf Since $O^p(L)\cap\Ker\,\varphi=\un$ by Corollary~\ref{condition BK}, the group $\Ker\,\varphi$ embeds into $L^{[p]}$, so it is a $p$-group. Let $F$ denote the Frattini subgroup of $\Ker\,\varphi$. Then $F$ is a normal subgroup of $L$. Moreover if $X$ is a subgroup of $L$ such that $XF=L$, then $F\leq\Ker\,\varphi\leq XF$, so $\Ker\,\varphi=(\Ker\,\varphi\cap X)F$, hence $\Ker\,\varphi\cap X=\Ker\,\varphi$, and then $XF=X=L$ since $F\leq\Ker\,\varphi\leq X$. It follows that $m_{L,F}=1$, thus $F=\un$ as $(L,\varphi)$ is a $B_K$-group. This shows that $\Ker\,\varphi$ is elementary abelian.\par
Let now $N=\mathop{\cap}_{P\in\mathcal{M}}\limits P$, where $\mathcal{M}$ is the set of normal subgroups of $L$ which are contained in $\Ker\,\varphi$, and maximal for these conditions (in other words the factor group $\Ker\,\varphi/P$ is a simple $\F_pL$-module). If $X$ is a subgroup of~$L$ such that $XN=L$, then $N\leq\Ker\,\varphi\leq XN$, so $\Ker\,\varphi=(\Ker\,\varphi\cap X)N$. But $\Ker\,\varphi\cap X$ is normalized by $X$ and $\Ker\,\varphi$, so it is normal in $L$. If $\Ker\,\varphi\cap X<\Ker\,\varphi$, then there is $P\in\mathcal{M}$ such that $\Ker\,\varphi\cap X\leq P$. Then $N\leq P$ also, and $\Ker\,\varphi=(\Ker\,\varphi\cap X)N\leq P$, contradicting $P<\Ker\,\varphi$. It follows that $m_{L,N}=1$, hence $N=\un$.\par
But then the product of the projection maps $\Ker\,\varphi\to\prod_{P\in\mathcal{M}}\limits\Ker\,\varphi/P$ is injective, and the latter is a semisimple $\F_pL$-module. Hence $\Ker\,\varphi$ is also a semisimple $\F_pL$-module. Now since $O^p(L)$ and $\Ker\,\varphi$ are normal subgroups of $L$ with trivial intersection, they centralize each other. In other words $\Ker\,\varphi$ is a module for the factor group $L^{[p]}=L/O^p(L)$. Then $\Ker\,\varphi$ is a semisimple $\F_pL^{[p]}$-module. As $L^{[p]}$ is a $p$-group, the action of $L^{[p]}$ on $\Ker\,\varphi$ has to be trivial. Hence $\Ker\,\varphi$ is central in $L$.\par
Let $Z$ be any subgroup of order $p$ of $\Ker\,\varphi$. Then $0=m_{L,Z}=1-\frac{k_L(Z)}{p}$, by Proposition 5.6.4 of~\cite{bisetfunctorsMSC}, where $k_L(Z)$ denotes the number of complements of $Z$ in $L$. It follows that $k_L(Z)=p$, so in particular there is a subgroup $H$ of $L$ such that $L=Z\times H$. Then the complements of $Z$ in $L$ are the groups of the form $\{\big(f(h),h\big)\mid h\in H\}$, where $f:H\to Z$ is any group homomorphism. It follows that there are exactly $p$ homomorphisms from $H$ to $Z\cong C_p$. Equivalently, there are exactly $p$ homomorphisms from the $p$-group $H^{[p]}$ to $C_p$, so $H^{[p]}$ is cyclic and non trivial. Since $\Ker\,\varphi$ embeds in $L^{[p]}\cong Z\times H^{[p]}$, the rank of $\Ker\,\varphi$ is at most~2.\par
We now observe that if $(M,\psi)\twoheadrightarrow (L,\varphi)$ is a surjective morphism of groups over $K$ - in particular if it is an isomorphism -, then $\psi(M)$ and $\varphi(L)$ are conjugate in $K$. Then there are three disjoint cases:
\begin{enumerate}
\item $\Ker\,\varphi=\un$. In this case, denoting by $\pi_H$ the inclusion map $H\hookrightarrow K$ and by $\varphi^0:L\to H$ the isomorphism induced by $\varphi$, we have $i\circ\varphi=\pi_H\circ\varphi^0$ for some inner automorphism $i$ of $K$ which conjugates $\varphi(L)$ to $H$. So $\varphi^0$ is an isomorphism from $(L,\varphi)$ to $(H,\pi_H)$ in $\grpover{K}$, and we are in Case 1 of Theorem~\ref{structure BKp groups}. 
\item $\Ker\,\varphi=Z\cong C_p$. Then we have seen that $L=Z\times H_1$, where $H_1$ is a subgroup of $L$ such that $H_1^{[p]}$ is cyclic and non trivial. In this case $\varphi$ induces an isomorphism $\varphi^0:H_1\to H=\varphi(L)$, and $\Id_Z\times \varphi^0$ is an isomorphism from $(L,\varphi)$ to $(Z\times H,j_H\circ\pi_H)$, where $\pi_H:Z\times H\to K$ is the projection onto $H$. Hence we are in Case~2 of Theorem~\ref{structure BKp groups}. 
\item $\Ker\,\varphi\cong C_p\times C_p$. Then let $Z$ be a subgroup of order $p$ of $\Ker\,\varphi$. Then we have seen that $L=Z\times H_1$, where $H_1$ is a subgroup of $L$ such that $H_1^{[p]}$ is cyclic and non trivial. In this case $Z_1=\Ker\,\varphi\cap H_1$ has order $p$, and $m_{L,Z_1}=0$ since $(L,\varphi)$ is a $B_K$-group. It follows that $Z_1$ must have also exactly $p$ complements in $L$. In particular, there is a subgroup $J$ of $L$ such that $L=Z_1\times J$. But then $Z_1\leq H_1\leq Z_1J$ implies that $H_1=Z_1\times H_2$, where $H_2=H_1\cap J$. Hence $L=Z\times Z_1\times H_2$, and moreover $H_2^{[p]}=\un$ since $H_1^{[p]}\cong Z_1\times H_2^{[p]}$ is cyclic. Then $\varphi$ induces an isomorphism $\varphi^0:H_2\to H=\varphi(L)$, and $\Id_{Z\times Z_1}\times \varphi^0$ is an isomorphism from $(L,\varphi)$ to $(Z\times Z_1\times H,j_H\circ\pi_H)$, where $\pi_H:Z\times Z_1\times H\to K$ is the projection onto $H$. Hence we are in Case 3 of Theorem~\ref{structure BKp groups}. 
\end{enumerate}
This completes the proof of Theorem~\ref{structure BKp groups}.\endpf
\begin{mth}{Corollary} Let $\sou{1}=\{0,1\}$ and $\sou{2}=\{0,1,2\}$ be totally ordered lattices of cardinality 2 and 3, respectively. Let $c_K$ (resp. $nc_K$) be the number of conjugacy classes of subgroups $H$ of $K$ such that $H^{[p]}$ is cyclic (resp. non cyclic). Then the lattice $\mathcal{I}_{\FB_K^{(p)}}$ of ideals of $\FB_K^{(p)}$ is isomorphic to the direct product of $c_K$ copies of $\sou{2}$ and $nc_K$ copies of $\sou{1}$. In particular it is a finite distributive lattice.
\end{mth}
\pf By Theorem~\ref{treillisp}, the lattice $\mathsf{e}_{L,\varphi}$ is distributive, isomorphic to the lattice $\mathcal{C}l_{\BgrKp{K}{p}}$ of closed subsets of $\BgrKp{K}{p}$. Moreover, the join-irreducible elements of $\mathcal{I}_{\FB_K^{(p)}}$ are the ideals $\mathsf{e}_{L,\varphi}$, for $(L,\varphi)\in \mathcal{I}_{\FB_K^{(p)}}$. By Theorem~~\ref{structure BKp groups}, the set $\BgrKp{K}{p}$ is finite, and contains three types of elements:
\begin{enumerate}
\item the elements $(H,j_H)$ of the first type, for $H\in[s_K]$.
\item the elements $(C_p\times H,j_H\circ\pi_H)$ of the second type, for $H\in[s_K]$ such that $H^{[p]}$ is cyclic and non-trivial.
\item the elements $(C_p\times C_p\times H,j_H\circ\pi_H)$ of the third type, for $H\in[s_K]$ such that $H^{[p]}$ is trivial.
\end{enumerate}
The only possible surjective morphisms between elements of $\BgrKp{K}{p}$ are of the following form:
\begin{itemize} 
\item $(C_p\times H,j_H\circ\pi_H)\twoheadrightarrow (H,j_H)$, where $H\in[s_K]$ is such that $H^{[p]}$ is cyclic and non trivial.
\item $(C_p\times C_p\times H,j_H\circ\pi_H)\twoheadrightarrow (H,j_H)$, where $H\in[s_K]$ is such that $H^{[p]}=\un$.
\end{itemize}
It follows that the poset $\BgrKp{K}{p}$ has as many connected components as conjugacy classes of subgroups of $K$. The connected components corresponding to subgroups $H$ for which $H^{[p]}$ is cyclic - trivial or not - are isomorphic to a totally ordered poset of size~2, and the other ones are posets with one element. Hence $\BgrKp{K}{p}$ is a disjoint union of $c_K$-components which are totally ordered of size~2, and $nc_K$ isolated points. The lattice of closed subsets of a totally ordered poset of size $n$ is a totally ordered lattice of size $n+1$, and the lattice of closed subsets of a disjoint union of posets is the direct product of the lattices of closed subsets of the pieces. This completes the proof.\endpf
\begin{rem}{Remark} \label{cite centrosp}As in Remark~\ref{cite centros}, it follows from Section 5.2.2 of~\cite{centros} that the category $\AMod{\FB_K^{(p)}}$ splits as a product
$$\AMod{\FB_K^{(p)}}\cong\prod_{H}\AMod{e_H^K\FB_K^{(p)}}\mvirg$$
of categories of modules over smaller Green biset functors $e_H^K\FB_K^{(p)}$, where $H\in[s_K]$. The above connected components correspond to this decomposition. In particular, when $H$ is a subgroup of $K$ such that $H^{[p]}$ is non cyclic, then the (commutative) Green functor $e_H^K\FB_K^{(p)}$ {\em has no non zero proper ideals}. It might therefore be called a {\em Green field}.
\end{rem}

\centerline{\rule{5ex}{.1ex}}
\begin{flushleft}
Serge Bouc - CNRS-LAMFA, Universit\'e de Picardie, 33 rue St Leu, 80039, Amiens Cedex 01 - France. \\
{\tt email : serge.bouc@u-picardie.fr}\\
{\tt web~~ : http://www.lamfa.u-picardie.fr/bouc/}
\end{flushleft}
\end{document}